\documentclass[article,3p]{elsarticle}

\usepackage{amsmath,amsfonts,amssymb,amsbsy,amscd}
\usepackage{mathrsfs,array}
\usepackage{verbatim}
\usepackage{amsthm}
\usepackage{graphicx}
\usepackage{tikz}
\usetikzlibrary{plotmarks}
\usetikzlibrary{arrows,shapes,trees}
\usepackage{subcaption}
\usepackage{hyperref}
 
\makeatletter
\def\ps@pprintTitle{%
 \let\@oddhead\@empty
 \let\@evenhead\@empty
 \def\@oddfoot{}%
 \let\@evenfoot\@oddfoot}
\makeatother

\newdefinition{rmk}{Remark}
\newproof{pf}{Proof}
\newproof{pot}{Proof of Theorem \ref{thm2}}

\def\B#1{\mbox{\boldmath{$#1$}}}

\newcommand{\pd}{\partial}

\newcommand{\bx}{\B{x}}

\newcommand{\onehalf}{\mbox{$\frac{1}{2}$}}
\newcommand{\VV}{\mathcal{V}}
\newcommand{\WW}{\mathcal{W}}

\def\be{\begin{equation}}
\def\ee{\end{equation}}
\def\ba{\begin{array}}
\def\ea{\end{array}}
\def\bea{\begin{eqnarray}}
\def\eea{\end{eqnarray}}
\def\beas{\begin{eqnarray*}}
\def\eeas{\end{eqnarray*}}
\newcommand{\bseq}{\begin{subequations}}
\newcommand{\eseq}{\end{subequations}}

\biboptions{numbers,sort&compress}

\usepackage{hyperref}
\hypersetup{
    colorlinks,
    citecolor=ck,
    filecolor=ck,
    linkcolor=ck,
    urlcolor=ck
}

\usepackage{changes}
\definechangesauthor[color=blue]{Rev.1}
\definechangesauthor[color=orange]{Rev.2}
\definechangesauthor[color=red]{Both Rev}
\definechangesauthor[color=purple]{Authors}

\journal{Computer Methods in Applied Mechanics and Engineering}

\begin{document}

\begin{frontmatter}

\title{
{\large{Correct energy evolution of stabilized formulations: The relation between VMS, SUPG and GLS via dynamic orthogonal small-scales and isogeometric analysis.}}\\ {\large{I: The convective--diffusive context
}}}
\author{M.F.P. ten Eikelder\corref{cor1}}
\cortext[cor1]{Corresponding author}
\ead{m.f.p.teneikelder@tudelft.nl}
\author{I. Akkerman}
\ead{i.akkerman@tudelft.nl}

\address{Delft University of Technology, Department of Mechanical, Maritime and Materials Engineering, P.O. Box 5, 2600 AA Delft, The Netherlands}
\date{today}

\begin{abstract}
This paper presents the construction of novel stabilized finite element methods in the convective--diffusive context that exhibit correct-energy behavior. Classical stabilized formulations can create unwanted artificial energy. Our contribution corrects this undesired property by employing the concepts of dynamic as well as  orthogonal small-scales within the variational multiscale framework (VMS). 
The desire for correct energy indicates that the large- and small-scales should be $H_0^1$-orthogonal.
Using this orthogonality the VMS method can be converted into the streamline-upwind Petrov-Galerkin (SUPG) or the Galerkin/least-squares (GLS) method.
Incorporating both large- and small-scales in the energy definition asks for dynamic behavior of the small-scales.
Therefore, the large- and small-scales are treated as separate equations.

Two consistent variational formulations which depict correct-energy behavior are proposed: (i) the Galerkin/least-squares method with dynamic small-scales (GLSD) and (ii) the dynamic orthogonal formulation (DO).
The methods are presented in combination with an energy-decaying generalized-$\alpha$ time-integrator. Numerical verification shows that dissipation due to the small-scales in classical stabilized methods can become negative, on both a local and a global scale. The results show that without loss of accuracy the correct-energy behavior can be recovered by the proposed methods. The computations employ NURBS-based isogeometric analysis for the spatial discretization.
\end{abstract}

\begin{keyword}
Correct-energy behavior \sep Stabilized methods \sep Projection \sep Orthogonal small-scales \sep Dynamic small-scales \sep Residual-based variational multiscale method \sep Convection-diffusion \sep Auxiliary flux \sep Isogeometric analysis \sep NURBS \sep Finite elements \sep Mixed formulation
\end{keyword}

\end{frontmatter}

\section{Introduction}
\label{sec:Introduction}

Stabilized methods and multiscale formulations form an auspicious, versatile and fundamental class of methodologies for finite element computations. The classical Galerkin variational formulation depicts correct-energy behavior although it has limitations concerning accuracy and stability. The popular stabilized methods, i.e. the Streamline upwind Petrov-Galerkin method (SUPG) \cite{BroHug82}, the Galerkin/least-squares method (GLS) \cite{HugFra89}, and the variational multiscale method (VMS) \cite{Hug95, Hug98}, overcome these issues, however show incorrect-energy behavior. In this paper we focus on convection-diffusion which serves as a model problem for more complex flow problems and turbulence.

This work is devoted to the construction of a new stabilized finite element method displaying \textit{correct-energy behavior}. Correct-energy behavior (or evolution) in a numerical method represents here that the method (i) does not create artificial energy and (ii) closely resembles the energy evolution of the continuous setting. A precise definition is included in Section \ref{sec:correct}. 

Our contribution fixes the incorrect energy deficiency by combining several ingredients. These are the \textit{dynamic} and \textit{orthogonal} behavior of the small-scales emerging from the stabilized methods, also referred to as \textit{dynamic orthogonal small-scales}, within the framework of isogeometric analysis.

\subsection{Dynamic small-scales}
In our quest for a correct-energy displaying formulation we learn that it is essential to use the so-called dynamic small-scales (also referred to as transient small-scales). This approach models the small-scales dynamically, i.e. with an ordinary differential equation in time, and takes its temporal contribution to the large-scale equation into account. This concept has originally been proposed by Codina in \cite{Cod02} and has been further analyzed in \cite{CodPriGuaBad07}.

\subsection{Orthogonal small-scales in VMS}
The multiscale stabilization method based on orthogonal small-scales serves as the next key ingredient of our approach. We link our choice of orthogonal small-scales to an optimality projector induced by the $H_0^1$-seminorm. This produces a highly attenuated and localized small-scale Green's function, which is very desirable property \cite{HugSan06}. We combine this methodology with residual-based variational multiscale modeling, a concept which emanates from VMS. The VMS approach finds many applications in incompressible turbulence, see e.g. \cite{StabRec2,HHOW04, BACHH07, BaCaCoHu07, LERC10, BaAk10}, and free surface flow \cite{ABKF11, LERC10}. Possible new directions in stabilized and multiscale methods are suggested in \cite{BaTaTez15}.

\subsection{Isogeometric analysis framework}
In addition, we employ the isogeometric analysis (IGA) methodology,
proposed by Hughes et al. in \cite{HuCoBa04}, which finds recent applications in various fields of science, see e.g. \cite{ABKF11, BCZH06, BHS12}. IGA is an effort to close the gap between on one hand
Computer-Aided Design (CAD) and on the other Computer-Aided Engineering (CAE). Finite element analysis (FEA) and CAD use a different representation for the geometry which makes a geometry update unpleasant and time-consuming. IGA corrects this deficiency by employing the same NURBS (non-uniform rational B-splines) geometry description as in CAD. This means that the NURBS surfaces in IGA match with the \textit{exact} CAD geometry, in contrast to FEA where the basis functions form an approximation of the CAD geometry. IGA leads to higher-order and higher-continuity discretizations on complex domains. Our choice for IGA is additionally motivated by the second derivatives in the weak formulations. This requires $C^1$-continuity of the basis functions. Furthermore, one of the main advantages of using the IGA methodology is that it guarantees the incompressibility constraint to hold exactly \cite{Evans13steadyNS, Evans13unsteadyNS}. This is a highly favorable property when the velocity field is not given, e.g. in case of the incompressible Navier-Stokes equations.

\subsection{Context}
The methodology is presented in the convective--diffusive model context which serves as a first step of this novel approach. The procedure is developed with the incompressible Navier-Stokes equations in mind which is the next step of this approach and is in itself presented in the sequel paper. In the context of stabilized methods a two-step approach, development for linear convection-diffusion followed by incompressible Navier-Stokes, is more common, see e.g. \cite{HugMal86, HugMal1986, Sha1091, FraFre92a, FraFre91}.

In the context of two-fluid flow phenomena, many numerical methodologies can unfortunately artificially create energy at the two-fluid interface. Even a small energy-inconsistency at the fluid surface can already lead to highly unstable behavior as is demonstrated in \cite{AkBaBeFaKe12}. To rectify this discrepancy, each of the components of the algorithm requires correct-energy behavior. When numerically solving air--water flow usually the components are (i) a standard incompressible Navier-Stokes solver and (ii) an algorithm describing the evolution of the air--water interface. Apart from its shared features with the incompressible Navier-Stokes equations, the convective--diffusive context is also required for the (level set) algorithm describing the evolution of the two-fluid interface.

\subsection{Outline}
The remainder of this paper is dedicated to the actual construction of a stabilized variational formulation which depicts correct-energy behavior and is summarized as follows. Section \ref{sec:continous} states the continuous form of the governing convection-diffusion equations, both in the strong form and in the weak form. The energy evolution linked to this formulation follows from the corresponding mixed-formulation which is derived with a Lagrange multiplier approach. Before proposing changes to existing stabilized methods, we introduce and discuss the energy evolution of the existing stabilized methods. Therefore, Section \ref{sec:stab} discusses the energy evolution in the standard VMS stabilized method with static small-scales. Section \ref{sec:correct} presents our quest towards a stabilized formulation depicting correct-energy behavior. It adds the two concepts (i) the dynamic behavior of the small-scales and (ii) the optimality projector yielding the vital orthogonality of the small-scales to the VMS formulation. Invoking the optimality projector in different ways leads to the other well-known stabilized formulations, namely SUPG and GLS. In Section \ref{sec:correct2} the demanded orthogonality between the small- and large-scales is enforced by the proper $H_0^1$-optimality projector. This yields an alternative variational multiscale stabilized formulation with correct-energy evolution. Furthermore, the methods demand a time-integrator which is correctly linked to an energy. Therefore, we re-address the generalized-$\alpha$ time-integration method. The energy demand leads to a certain parameter family of the time-stepping parameters. Section \ref{sec:time} discusses this temporal-integration method. Section \ref{sec:Numerical verification} presents numerical verification while employing NURBS basis functions for the computations. In Section \ref{sec:conc}, we draw conclusions and outline avenues for future research.

\section{The continuous convection-diffusion equation}\label{sec:continous}

\subsection{Strong formulation}\label{sec:strong}

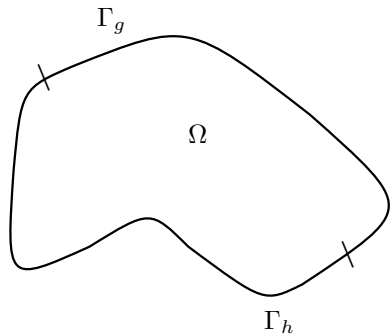
\begin{figure}[h!]
\centering
\begin{tikzpicture}
\draw[line width=0.3mm, black ] (2,0.4) .. controls (3.5,1) .. (5.1, -0.25) .. controls (6.5,-1.5) .. (5, -2.5) .. controls (4.5,-2.75) .. (3.5, -2).. controls (3.0,-1.5) .. (2.2, -2) .. controls (1.1,-2.5) .. (1.2, -1.0) .. controls (1.3,0.1) .. (2.0,0.4);
\node[text width=3cm] at (5.0,-0.5) {$\Omega$};
\node[text width=3cm] at (6.0,-3.0) {$\Gamma_h$};
\node[text width=3cm] at (3.8,1.0) {$\Gamma_g$};
\node[text width=3cm] at (7.0,-2.1) {$\boldsymbol{\backslash}$};
\node[text width=3cm] at (3.0,0.25) {$\boldsymbol{\backslash}$};
\end{tikzpicture}
  \caption{Spatial domain $\Omega$ with its boundaries $\Gamma = \Gamma_g \cup \Gamma_h$.}
  \label{fig:domain1}
\end{figure}

Let $\Omega$ denote the spatial domain with boundary $\Gamma = \Gamma_g \cup \Gamma_h $, see Figure \ref{fig:domain1}.
\\

The governing equations of the convection-diffusion problem in strong form read
\begin{subequations}
  \label{sec:GE, subsec:SF, CD Strong}
  \begin{alignat}{1}
    \pd_t \phi + \B{a} \cdot\nabla \phi -\nabla \cdot \kappa\nabla \phi = f & \quad \text{in} \quad \Omega\times \mathcal{I}, \label{sec:GE, subsec:SF, CD Strong, CD eq}\\
    \phi = g & \quad \text{in} \quad \Gamma_g \times \mathcal{I}, \label{sec:GE, subsec:SF, CD Strong, Dir eq}\\
 - a_n^- \phi   +\kappa  \pd_n \phi  = h & \quad \text{in} \quad \Gamma_{h}  \times \mathcal{I}, \label{sec:GE, subsec:SF, CD Strong, Neu eq}\\
    \phi(\bx,0) = \phi_0(\bx)  & \quad \text{in} \quad \Omega, \label{sec:GE, subsec:SF, CD Strong, IC}
  \end{alignat}
\end{subequations}
where $t \in \mathcal{I}=(0,T)$ is the time with final time $T>0$, $x \in \Omega$ the spatial coordinate, $\phi=\phi(\bx,t): \Omega \times \mathcal{I}$ the dependent variable with time derivative $\pd_t \phi$, normal flux $\pd_n \phi = \B{n} \cdot\nabla \phi$ and $f: \Omega \times \mathcal{I} \rightarrow \mathbb{R}$, $g: \Gamma_g \times \mathcal{I} \rightarrow \mathbb{R}$, $h: \Gamma_h \times \mathcal{I} \rightarrow \mathbb{R}$  and $ \phi_0: \Omega \rightarrow \mathbb{R}$ are prescribed data. The convective velocity $\B{a}=\B{a}(\bx)$ is a given solenoidal vector field ($\nabla \cdot \B{a} = 0$) and $\kappa \geq 0$ denotes the given diffusivity. The outward unit normal to $\Gamma$ is $\B{n}$ and the normal velocity component denotes $a_n=\B{a}\cdot\B{n}$ with positive and negative parts $a_n^{\pm}=\tfrac{1}{2}(a_n\pm|a_n|)$.

\subsection{Weak formulation}\label{sec:weak}

Let $\WW^0$ and $\WW^g$ denote suitable function-spaces
satisfying the homogeneous and non-homogeneous Dirichlet conditions, respectively.
Using these spaces the standard weak formulation of the problem reads as follows:\\

\textit{Find $\phi \in \WW^g$ such that for all $w \in \WW^0$,}
\begin{align}
  \label{sec:GE, subsec:SWF, standard weak form}
 \left(w,\pd_t \phi\right)_{\Omega} +  \left( w, \B{a} \cdot\nabla \phi\right)_{\Omega}
 - \left(w, a_n^- \phi \right)_{\Gamma_{h}} + \left(\nabla w,\kappa\nabla \phi \right)_{\Omega}=(w,f)_{\Omega}+(w,h)_{\Gamma_h}.
\end{align}
Here $\left(\cdot,\cdot\right)_D$ denotes the $L^2(D)$ inner product over $D$. 
Consistency of the strong (\ref{sec:GE, subsec:SF, CD Strong}) and the weak formulation (\ref{sec:GE, subsec:SWF, standard weak form}) easily follows from applying integration by parts on the diffusive term.

Instead of enforcing the Dirichlet boundary conditions \textit {a priori}, it is also possible to relax this condition in the function space by employing a Lagrange multiplier setting.
The weak statement translates into a \textit{mixed formulation}: \\

\textit{Find $\left(\phi,\lambda_\Omega \right) \in \WW \times \VV$ such that for all $\left(w,q\right) \in \WW\times \VV$,}
\begin{subequations}
  \label{sec:GE, subsec:SWF, LM weak form}
  \begin{alignat}{1}
 (w,\lambda_\Omega )_{\Gamma_g}=&\left(w,\pd_t \phi\right)_{\Omega}
 + \left(w,\B{a} \cdot\nabla \phi\right)_{\Omega}
 -\left(w, a_n^- \phi \right)_{\Gamma_{h}}+ \left(\nabla w,\kappa\nabla \phi \right)_{\Omega}
 -(w,f)_{\Omega} - (w, h )_{\Gamma_h}, \\
 (q, \phi )_{\Gamma_g}  =& (q,  g )_{\Gamma_g}.
\end{alignat}
\end{subequations}
Here $\WW$ represents the unrestricted function space and $\VV$ is a suitable Lagrange multiplier space. Consult \cite{HEML00, HuWe05} for the appropriate construction of the spaces. The following Section employs this formulation to derive energy statements.

Applying an appropriate choice of weighting functions $w$ and $q$ in (\ref{sec:GE, subsec:SWF, LM weak form}) and subsequently performing a partial integration step recovers the strong form (\ref{sec:GE, subsec:SF, CD Strong}). Additionally, the expression for the Lagrange multiplier follows as a complimentary result
\begin{align}
 \lambda_\Omega = \kappa \pd_n \phi,
\end{align}
and equals the diffusive flux. Note that the continuous setting allows us to provide a closed form of the Lagrange multiplier. This does not hold in a discrete setting. Furthermore, the subscript in the notation of the Lagrange multiplier is added for consistency with Section \ref{sec:loc_energy}.

\subsection{Global energy evolution}\label{sec:glob_energy}

The evolution of the energy linked to the strong form (\ref{sec:GE, subsec:SF, CD Strong}) follows from choosing $w=\phi$ and $q=\lambda$ in the mixed formulation (\ref{sec:GE, subsec:SWF, LM weak form}). Subtracting the resulting equations yields
\begin{align}
  \label{eq:glb_en_0}
  \left(\phi,\pd_t \phi\right)_{\Omega}   + \left(\nabla \phi,\kappa\nabla \phi \right)_{\Omega}
  + \left(\phi,\B{a} \cdot\nabla \phi\right)_{\Omega}
  - \left(\phi, a_n^-\phi \right)_{\Gamma_{h}}
    = (g,\lambda_\Omega )_{\Gamma_g}
  +(\phi,f)_{\Omega}  +(\phi,h)_{\Gamma_h}.
\end{align}
By performing integration by parts on the interior convective term and employing the divergence-free constraint, the convective term turns into a boundary term. The energy, defined as $E_{\Omega}=\frac{1}{2}(\phi,\phi)_{\Omega}$, evolves as
\begin{align}
  \label{eq:glb_en}
  \dfrac{{\rm d}}{{\rm d}t} E_{\Omega}=-\|\kappa^{1/2}\nabla \phi\|_{\Omega}^2 +(\phi,f)_{\Omega}
  - (1,F_\Omega)_{\Gamma},
\end{align}
where $||\cdot||_{D}$ defines the standard $L^2$-norm over $D$. The conservative energy flux provides a different contribution on each segment of the boundary:
\begin{align}
  \label{eq:glb_flux}
F_\Omega =
\left \{ \begin{array}{cc}
  a_ne  -  g \lambda_\Omega   &\text{on}\quad \Gamma_g,  \\
  |a_n|e  - \phi h&\text{on} \quad \Gamma_h,\\
  0 & \quad\text{elsewhere},
  \end{array} \right.
\end{align}
with $e:=\onehalf \phi^2$ the pointwise energy. The terms on the Dirichlet boundary are (i) the amount of energy flowing in and out by convection and (ii) the energy gained and lost by diffusion through the boundary, respectively. On the other boundary, the terms represent (i) the energy loss by the strict convective outflow and (ii) the energy change by the flux boundary condition. The energy can only increase as a result of the prescribed body force or the boundary conditions, represented by the last two terms on the right-hand side of (\ref{eq:glb_en}). The diffusive term, when active, contributes to a decay of the energy. The last term on the right-hand side of (\ref{eq:glb_en}) represents the convective and diffusive fluxes of energy across the interface. Substitution of the boundary condition and the Lagrange multiplier (again possible because of the continuous setting) and a partial integration step leads to the alternative expression of the flux
\begin{align}
  \label{eq:glb_flux}
F_\Omega =  a_ne  -  \kappa \pd_n e.
\end{align}
The two terms respectively describe the convective and viscous-driven flow of energy.
\subsection{Localized energy evolution}\label{sec:loc_energy}

\begin{figure}[h!]
\centering
\begin{tikzpicture}
\draw[line width=0.3mm, black ] (2,0.4) .. controls (3.5,1) .. (5.1, -0.25) .. controls (6.5,-1.5) .. (5, -2.5) .. controls (4.5,-2.75) .. (3.5, -2).. controls (3.0,-1.5) .. (2.2, -2) .. controls (1.1,-2.5) .. (1.2, -1.0) .. controls (1.3,0.1) .. (2.0,0.4);
\draw[line width=0.3mm, black ] (5.1, -0.25) .. controls (3.6,-0.2) .. (3.5, -2);
\node[text width=3cm] at (3.5,-0.5) {$\Omega-\omega$};
\node[text width=3cm] at (6.0,-1.5) {$\omega$};
\node[text width=3cm] at (6.0,-3.0) {$\Gamma_h(\omega)$};
\node[text width=3cm] at (7.7,-1.0) {$\Gamma_g(\omega)$};
\node[text width=3cm] at (3.8,1.0) {$\Gamma_g(\Omega-\omega)$};
\node[text width=3cm] at (2.6,-2.8) {$\Gamma_h(\Omega-\omega)$};
\node[text width=3cm] at (5.3,-0.1) {$\chi_\omega$};
\node at (5cm,-1) {\pgfuseplotmark{cross*}};
\node[text width=3cm] at (7.0,-2.1) {$\boldsymbol{\backslash}$};
\node[text width=3cm] at (3.0,0.25) {$\boldsymbol{\backslash}$};
\node[text width=3cm] at (6.5, -0.25) {$\boldsymbol{/}$};
\node[text width=3cm] at (4.9, -2) {$\boldsymbol{/}$};
\end{tikzpicture}
  \caption{Spatial domain $\Omega$ with a subdomain $\omega\subset \Omega$. The shared boundary of $\omega$ and its complement is $\chi_\omega$. The boundaries $\Gamma_g$ and $\Gamma_h$ split according to $\omega$.}
  \label{fig:domain2}
\end{figure}
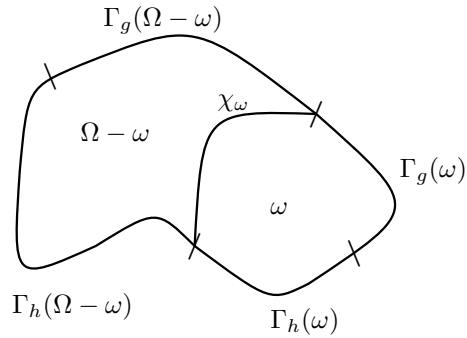

This Section presents a localized version of (\ref{eq:glb_en}) for arbitrary subdomains $\omega \subset \Omega$ with boundary $\pd \omega$. The complement domain denotes $\Omega - \omega$ with boundary $\pd (\Omega-\omega)$ and the shared boundary of the both subdomains is $\chi_\omega := \pd \omega \cap \pd (\Omega-\omega)$. Figure \ref{fig:domain2} shows the domain with its boundaries.

The variational statement consists of the weak formulation (\ref{sec:GE, subsec:SWF, LM weak form}) enforced on the subdomain $\omega$ and is again augmented with a Lagrange multiplier that enforces the Dirichlet boundary condition. The unrestricted solution space $\WW_\omega$ allows discontinuities across the subdomain interface and the space of suitable Lagrange multipliers denotes $\VV_\omega$. The weak statement reads:\\

\textit{Find $\left(\phi, \lambda_{\omega} \right) \in  \WW_\omega \times \VV_\omega$ such that for all $\left(w, q\right) \in \WW_\omega \times \VV_\omega$,}
\begin{subequations}
 \label{eq:loc_weak}
  \begin{alignat}{2}
  (w, \lambda_{\omega})_{\chi_\omega} +(w,\lambda_{\omega} )_{\Gamma_{g}(\omega)} =&\left(w,\pd_t \phi\right)_{\omega}
  + \left(\nabla w,\kappa\nabla \phi \right)_{\omega}
 + \left(w,\B{a} \cdot\nabla \phi\right)_{\omega}\nonumber \\[10pt]
 &- \left(w, a_n^- \phi\right)_{\Gamma_{h}(\omega)}
   -   (w,f)_{\omega}
    - (w, h )_{\Gamma_{h}(\omega)},\\[10pt]
      (q^h,[\phi^h])_{\chi_\omega} + (q, \phi )_{\Gamma_g(\omega)}=& (q,  g )_{\Gamma_g(\omega)} ,
  \end{alignat}
\end{subequations}
where $\Gamma_g(D):=\Gamma_g\cap \partial D$ and $\Gamma_h(D):=\Gamma_h\cap \partial D$ for domain $D$. The jump term $[\phi^h]$ is given by
\begin{equation}
 [\phi^h]:=\phi^h_\omega -\phi^h_{\Omega-\omega},
\end{equation}
where the terms are defined on $\omega$ and $\Omega-\omega$, respectively.
The Lagrange multiplier can be interpreted as an auxiliary flux across the interface $\chi_\omega$, it represents the flow outward $\omega$. The weak form of the complement domain easily follows by replacing $\omega$ by $\Omega-\omega$ in (\ref{eq:loc_weak}). The superposition of the both formulations leads to the balance
\begin{equation}
  \lambda_{\omega} + \lambda_{\Omega-\omega} = 0.
\end{equation}
Thus that what flows out of $\omega$ enters in $\Omega-\omega$. See \cite{HEML00} for the formal details of such a derivation. Again, a partial integration step provides the expression for the Lagrange multipliers:
\begin{subequations}
  \label{eq:lm ex}
  \begin{alignat}{2}
    \lambda_{\omega} =& \kappa \pd_{n_\omega} \phi,\\
    \lambda_{\Omega-\omega} =& \kappa \pd_{n_{\Omega-\omega}} \phi,
  \end{alignat}
\end{subequations}
with $\pd_{n_D}$ the directional derivative outward of a domain $D$.
The local energy statement follows when choosing $w=\phi$ and $q=\lambda_{\omega}$, this yields
\begin{align}
  \label{eq:loc_en_0}
   (\phi, \lambda_{\omega})_{\chi_\omega}+ (g,\lambda_{\omega} )_{\Gamma_g(\omega)}=&  \left(\phi,\pd_t \phi\right)_{\omega}
+ \|\kappa^{1/2}\nabla \phi\|_{\omega}^2
 + \left(\phi,\B{a} \cdot\nabla \phi\right)_{\omega}\nonumber \\[8pt]
    &- \left(w, a_n^- \phi\right)_{\Gamma_h(\omega)}
  -   (\phi,f)_{\omega}        - (\phi, h )_{\Gamma_h(\omega)}.
 \end{align}
By applying integration by parts on the convective term we find that the energy on subdomain $\omega$ evolves as
\begin{align}
 \label{eq:loc_en}
  \dfrac{{\rm d}}{{\rm d}t} E_{\omega}=- \|\kappa^{1/2}\nabla \phi\|_{\omega}^2+(\phi,f)_{\omega}
  - (1,  F_\omega)_{\partial \omega},
\end{align}
where the outward energy flux now splits into three parts
\begin{align}
 \label{eq:loc_flux}
F_\omega =
\left \{ \begin{array}{ccl}
  a_{n_\omega}e  - g \lambda_\omega        &&\text{on}\quad \Gamma_g(\omega),  \\
  |a_{n_\omega}|e  - \phi h              &&\text{on} \quad \Gamma_h(\omega), \\
   a_{n_\omega}e   - \phi \lambda_\omega   &&\text{on}\quad\chi_\omega,\\
  0 &&\text{elsewhere}.
  \end{array} \right.
\end{align}
In comparison with global energy behavior, the additional last term represents an energy flux, with a convective and diffusive component, across the subdomain interface $\chi_\omega$. Again, the substitution of the boundary condition and the Lagrange multiplier, and subsequently performing a partial integration step lead to
\begin{align}
  \label{eq:glb_flux}
F_\omega =  a_{n_\omega}e  -  \kappa \pd_{n_\omega} e \quad \text{on}\quad \partial \omega.
\end{align}

\subsection*{Remark}
This Section provides all the statements in a continuous form. A direct consequence is that the standard discrete setting, the Galerkin method, displays correct-energy behavior.\\

This paper now presents the energy evolution of standard stabilized methods and subsequently constructs a methodology that closely resembles the local energy evolution of the continuous equation. In particular, the design of the method precludes artificial local energy creation.

\subsection*{Remark}
To increase the readability of this paper we now restrict ourselves to boundary conditions precluding the energy flux $F$ on $\Gamma$. This occurs for example when employing homogeneous Dirichlet and periodic boundary conditions. The proposed methodology can easily be generalized to domains with non-homogeneous boundaries.

\section{Energy evolution of the variational multiscale approach}\label{sec:stab}

This Section concerns the energy evolution in the stabilized residual-based variational multiscale (RBVMS) formulation. Therefore we start off with a brief recap of the canonical VMS method.

\subsection{The multiscale split}\label{sec:vms1}
The residual-based variational multiscale approach emanates, as the name suggests, from the theory of the variational multiscale methods. This approach explicitly treats the solution component not be represented by the discretization in an approximate sense. This is done as follows. The trial solution and weighting function spaces split as
\begin{align}
  \WW = \WW^h \oplus \WW',
\end{align}
where $\WW^h$ is the space spanned by the finite-dimensional discretization and $\WW'$ is its infinite-dimensional complement in $\WW$. Based on the multiscale split in the space the components of the solution and weight decouple as
\begin{align}
  \phi =& \phi^h + \phi',\nonumber\\
  w =& w^h + w',
\end{align}
with $\phi^h, w^h \in \WW^h$ and $\phi', w' \in \WW'$ the large-scales and the small-scales solution and weight, respectively. The small-scale space $\WW'$ represents the component of $\WW$ not reproduced by the grid and is therefore also called space of fine-scales, sub-scales or subgrid-scales\footnote{The terms sub-scales or subgrid-scales could be linked to a turbulence modeling character of the approach. The current method does not fit in that framework. To emphasize this difference we use the terminology small-scales.}. In order to obtain a well-defined decomposition for a given $v \in \WW$, the elements $v^h \in \WW^h, v' \in \WW'$ with $v=v^h+v'$ require a unique definition. Employing an optimality projector $\mathscr{P}^h: \WW \rightarrow \WW^h$ for the decoupling as\footnote{There are infinitely many choices for the projector $\mathscr{P}^h$. Linear projectors suffice for the current problem. More details can be found in Hughes \cite{Hug95}.}:
\begin{align}
\label{eq:opt_proj}
  v^h =& \mathscr{P}^h v,\nonumber\\
  v' =& \left(\mathscr{I}-\mathscr{P}^h\right) v,
\end{align}
achieves uniqueness. Here $\mathscr{I}: \WW \rightarrow \WW$ is the identity operator. Using this multiscale split we arrive at the following alternative -- equivalent -- weak statement:\\

\textit{Find $\phi^h \in \mathcal{W}^h, \phi' \in \mathcal{W}'$ such that for all $w^h \in \WW^h, w' \in \WW'$,}
\begin{subequations}
  \label{eq:vms_ex}
  \begin{alignat}{1}
   (w^h,\pd_t \phi^h+\B{a} \cdot\nabla \phi^h)_{\Omega}  + (\nabla w^h, \kappa\nabla \phi^h)_{\Omega}
  + (w^h, \pd_t \phi'+\B{a} \cdot\nabla \phi')_{\Omega} + (\nabla w^h, \kappa\nabla \phi')_{\Omega}
  =& (w^h, f)_{\Omega},   \label{eq:large} \\
    (w',\pd_t \phi^h + \B{a} \cdot\nabla \phi^h)_{\Omega}  + (\nabla w', \kappa\nabla \phi^h)_{\Omega}
  + (w',\pd_t \phi'+\B{a} \cdot\nabla \phi')_{\Omega} + (\nabla w', \kappa\nabla \phi')_{\Omega} =& (w', f)_{\Omega}. \label{eq:small}
  \end{alignat}
\end{subequations}
Note that this formulation is still exact. However, the space  $\WW'$ is infinite-dimensional and as such not amenable for a discrete implementation.

\subsection{The VMS numerical formulation}\label{sec:vms2}

The weak formulation (\ref{eq:vms_ex}) converts into a numerical formulation when the small-scale equation (\ref{eq:small}) is replaced by an approximation for the small-scale solution $\phi'$. The small-scale equation can be written in the form
\begin{align}\label{eq:small2}
  (w',\pd_t \phi'+\B{a} \cdot\nabla \phi' - \kappa \Delta \phi')_{\Omega}=-(w',\mathscr{R}\phi^h)_{\Omega},
\end{align}
where the large-scale residual reads
\begin{align}
  \mathscr{R}\phi^h=\pd_t \phi^h+\B{a} \cdot\nabla \phi^h - \kappa \Delta \phi^h-f.
\end{align}
This implies that the small-scales are driven by the residuals of the large-scales. The corresponding Euler-Lagrange form of the small-scale equation reads
\begin{align}\label{eq:EL}
  \pd_t \phi'+\B{a} \cdot\nabla \phi' - \kappa \Delta \phi'=&  - \mathscr{R}\phi^h.
\end{align}
Note that the Euler-Lagrange equations are in strong form, i.e. the weight $w'$ in (\ref{eq:small2}) is ignored. This pertains to both the small-scale solution as well as the residual forcing.

Employing a Green's function provides an explicit expression for the small-scales. In this expression the integral operator is approximated by an algebraic stabilization parameter $\tau_{\text{static}}$. This step is necessary to arrive at an implementable method. Details of this approximation can be found in \cite{Hug95}. To obtain a stabilized formulation, the small-scales are modeled as:
\begin{subequations}
  \label{eq:small_stat}
  \begin{alignat}{1}
  \hat{\phi}'= &-\tau_{\text{stat}} \mathscr{R}\phi^h,\label{eq:small_stat, algebraic}\\
    \pd_t \hat{\phi}' =& 0, \label{eq:small_stat, der}
\end{alignat}
\end{subequations}
where $\tau_{\text{stat}}$ is a positive stabilization parameter. In the following we ignore the hat-sign. This algebraic operator depends on both the
physics and the discretization. More details can be found in Hughes and Sangalli \cite{HugSan06}.

The definition of the stabilization parameter is inspired by the theory of stabilized methods for convection-diffusion equations (see e.g. \cite{HugMal86, Sha1091}), and
reads:
\begin{align} \label{eq:tau_static}
  \tau_{\text{stat}}=\left(\tau^{-2}_{\text{conv}}+\tau^{-2}_{\text{diff}} + \tau^{-2}_{\text{time}}\right)^{-1/2},
\end{align}
where
\begin{subequations}
  \begin{alignat}{1}
  \tau_{\text{conv}}^{-2}=& \mathbf{a}\cdot \mathbf{G} \mathbf{a},\\
  \tau_{\text{diff}}^{-2}=& C_I\kappa^2 \mathbf{G}:\mathbf{G},\\
  \tau_{\text{time}}^{-2}=& \left(\frac{\alpha_m}{\alpha_f \gamma \Delta t}\right)^2.
\end{alignat}
\end{subequations}
Here $\mathbf{G}$ is the second-rank metric tensor given by
\begin{align}
  \mathbf{G}=\frac{\pd \boldsymbol{\xi}}{\pd \bx}^T\frac{\pd \boldsymbol{\xi}}{\pd \bx},
\end{align}
where $\pd \boldsymbol{\xi}/\pd \bx$ is the inverse Jacobian of the map between the elements in the reference and physical domain.
The stabilization parameter treats deformed and curved domains naturally due to its direct dependence on $\mathbf{G}$. The metric tensor $\mathbf{G}$ scales as $h^{-2}$ where $h$ is the mesh size. The positive constant $C_I$
is defined by an inverse estimate. It is independent of the mesh size and can be computed from
an element-wise eigenvalue problem  \cite{HarHug92}.

The definition of $\tau_{\text{time}}$ is based on the generalized-$\alpha$ time-integrator given in Section \ref{sec:time}. Here $\alpha_f, \alpha_m, \gamma$ are algorithmic time-stepping
coefficients and $\Delta t$ is the time step.
It reduces to the commonly used value of $4/\Delta t^2$ when employing $\rho_\infty =1$, see e.g. \cite{BACHH07, BaAk10, BaCaCoHu07} (consult these references for the definition of $\rho_\infty$). This choice results in the Crank-Nicolson time-integrator, see Section \ref{sec:time, subsec:TIFSE}.

Employing integration by parts on the stabilized terms, the small-scales appear without derivatives. The resulting form is the VMS stabilized statement\\

\textit{Find $\phi^h \in \mathcal{W}^h$ such that for all $w^h \in \WW^h$,}
\begin{subequations}
\label{eq:VMS_weak}
  \begin{alignat}{1}
  (w^h,\pd_t \phi^h )_{\Omega}  + (w^h,\B{a} \cdot\nabla \phi^h)_{\Omega} + (\nabla w^h, \kappa\nabla \phi^h)_{\Omega}
   - (\B{a} \cdot\nabla w^h + \kappa \Delta w^h,  \phi')_{\tilde{\Omega}}
  =& (w^h, f)_{\Omega},   \\
  \tau_{\text{stat}}^{-1} \phi'=& -\mathscr{R}\phi^h. \label{VMS stat ss}
\end{alignat}
\end{subequations}
Here we have subdivided the domain $\Omega$ into elements $\Omega_e$. The domain of element interiors $\tilde{\Omega}$ does not include the element boundaries and reads
\begin{align}
  \tilde{\Omega} = \displaystyle \bigcup_e \Omega_e.
\end{align}
It is important to emphasize to that we treat the small-scale expression (\ref{VMS stat ss}) as a separate equation. At this stage a straightforward substitution is certainly possible, however when the small-scales are modeled dynamically, this is not the case anymore. In line with the analysis in later sections we therefore omit substitution here.
\subsection{Local energy evolution of the VMS formulation}\label{sec:stab, subsec: Local energy evolution of the VMS formulation}

To arrive at local energy evolution, we augment the weak formulation in Lagrange multiplier setting form to allow discontinuous functions across subdomains, similar to (\ref{eq:loc_weak}). The weak statement reads: \\

\textit{Find $\left( \phi^h, \lambda^h_\omega \right) \in \WW^h_\omega \times \VV^h_\omega$ such that for all $\left( w^h, q^h\right) \in \WW^h_\omega \times \VV^h_\omega$,}
\begin{subequations}
\label{sec:VMS_weak_LM}
  \begin{alignat}{1}
  (w^h, \lambda^h_\omega)_{\chi_\omega} =& (w^h,\pd_t \phi^h )_{\omega}  + (w^h,\B{a} \cdot\nabla \phi^h)_{\omega} + (\nabla w^h, \kappa\nabla \phi^h)_{\omega}
   - (w^h,f)_{\omega}\nonumber \\&
   - (\B{a} \cdot\nabla w^h + \kappa \Delta w^h,  \phi')_{\tilde{\omega}}
  ,\label{sec:VMS_weak_LM, large scale}\\
  (q^h,[\phi^h])_{\chi_\omega}=& 0,\\
  \tau_{\text{stat}}^{-1} \phi'=& -\mathscr{R}\phi^h.
\end{alignat}
\end{subequations}
Here $\tilde{\omega}$ represents the domain of element interiors of $\omega$. The discretization does not allow explicit evaluation of the Lagrange multiplier.
We select $w^h =\phi^h$ in the large-scale equation (\ref{sec:VMS_weak_LM, large scale}) and add the small-scale equation multiplied by $\phi'$ and integrate. The resulting statement is:
\begin{align}\label{eq:VMS_energy_0}
  (\phi^h ,\pd_t \phi^h)_{\omega}
  +(\phi',\pd_t \phi^h)_{\tilde{\omega}}
  + \| \tau_{\text{stat}}^{-1/2} \phi'\|^2_{\tilde{\omega}}
  +  \|\kappa^{1/2}\nabla \phi^h \|^2 _{\omega}
  - (\phi^h ,f)_{\omega} & \nonumber \\
  -(\phi',f)_{\tilde{\omega}}-2(\kappa\Delta \phi^h,\phi')_{\tilde{\omega}}
   +\onehalf (\phi^h , a_n \phi^h)_{\chi_\omega}
  -(\phi^h, \lambda^h_\omega)_{\chi_\omega} &= 0.
\end{align}
Here we have employed the incompressibility constraint to convert the interior convective term to a boundary term.
\subsection*{Remark}
When the velocity field is obtained by a numerical method the incompressibility constraint is often not exactly fulfilled though. However, by appropriately employing isogeometric analysis this can be achieved exactly \cite{Evans13unsteadyNS}. Our implementation already employs the proper IGA spaces to allow a smooth transition to the incompressible Navier-Stokes equations.\\

The \textit{local large-scale energy} is the energy of resolved solution: $E^h_\omega=\frac{1}{2}\left(\phi^h,\phi^h\right)_{\omega}$ and evolves by (\ref{eq:VMS_energy_0}) as:
\begin{align}\label{eq:VMS_energy resolved}
 \dfrac{{\rm d}}{{\rm d}t} E^h_{\omega}=&  - \| \kappa^{1/2} \nabla \phi^h \|^2_{\omega} + (\phi^h ,f)_{\omega} -(1,F^h_\omega)_{\chi_\omega}
    \nonumber  \\
 & - \| \tau_{\text{stat}}^{-1/2}\phi'\|^2_{\tilde{\omega}} + (\phi' ,f)_{\tilde{\omega}} + 2(\kappa\Delta \phi^h,\phi')_{\tilde{\omega}}-(\phi',\pd_t \phi^h)_{\tilde{\omega}},
 \end{align}
with the energy flux
\begin{equation}\label{energy flux VMS}
  F_\omega^h = a_ne^h - \lambda^h_\omega\phi^h.
\end{equation}
where the \textit{pointwise large-scale energy} is $e^h:=\onehalf \phi^h\phi^h$. The \textit{local total energy} is defined using the superposition of the small-scales and large-scales as: $E_\omega=\frac{1}{2}\left(\phi^h+\phi',\phi^h+\phi' \right)_{\tilde{\omega}}$. Its evolution directly follows:
\begin{align}\label{eq:VMS_energy}
 \dfrac{{\rm d}}{{\rm d}t} E_{\omega}=&  - \| \kappa^{1/2} \nabla \phi^h \|^2_{\omega} + (\phi^h ,f)_{\omega} -(1,F^h_\omega)_{\chi_\omega}
    \nonumber  \\
 & - \| \tau_{\text{stat}}^{-1/2}\phi'\|^2_{\tilde{\omega}} + (\phi' ,f)_{\tilde{\omega}} + 2(\kappa\Delta \phi^h,\phi')_{\tilde{\omega}}+  (\pd_t \phi',\phi^h+\phi')_{\tilde{\omega}}.
 \end{align}
We observe from this relation that the standard static VMS formulation does not possess a desirable energy behavior. The first line closely resembles the continuous energy evolution relation. No explicit expression for $\lambda^h_{\omega}$ exists. The second line appears as a result of the stabilization terms. Its first term contributes to a decay of the energy, which is favorable from a stability argument. It can be interpreted as the diffusive energy decay of the missing small-scales. The last two terms are problematic. These unsymmetric terms can be bounded by both the physical diffusion $\| \kappa^{1/2} \nabla \phi^h \|^2_{\omega}$ and the numerical diffusion $\| \tau_{\text{stat}}^{-1/2}\phi'\|^2_{\tilde{\omega}}$. The procedure is analogous to the standard coercivity analysis: apply Cauchy-Schwarz and Young's inequality subsequently. This leads to restrictions on the stabilization parameter $\tau_{\text{stat}}$ depending on the diffusivity and the time step. More importantly, the overall diffusion of the method can be less than the physical diffusion. Hence, the small-scales can artificially create energy, which we numerically show in Section \ref{sec:Numerical verification}, and are therefore both numerically and physically undesirable. The next section corrects this deficiency. 
\subsubsection*{Remark}
The global energy evolution easily follows when substituting $\omega = \Omega$ and $\tilde{\omega} = \tilde{\Omega}$ into (\ref{eq:VMS_energy}):
\begin{align}
  \label{eq:VMS global energy}
 \dfrac{{\rm d}}{{\rm d}t} E_{\Omega}=&  - \| \kappa^{1/2} \nabla \phi^h \|^2_{\Omega} + (\phi^h ,f)_{\Omega}
    \nonumber  \\
 & - \| \tau_{\text{stat}}^{-1/2}\phi'\|^2_{\tilde{\Omega}} + (\phi' ,f)_{\tilde{\Omega}}  + 2(\kappa\Delta \phi^h,\phi')_{\tilde{\Omega}}+  (\pd_t \phi',\phi^h+\phi')_{\tilde{\Omega}}.
\end{align}
Note the cancellation of the local energy flux.

\section{Toward a stabilized formulation with correct-energy behavior}
\label{sec:correct}
This Section presents a path with alternative stabilized formulations towards rectification of the discrepancy indicated in the previous section. First we adopt the concept of dynamic small-scales to eliminate the unwanted terms containing the temporal derivatives. Next, the undesirable diffusive term vanishes when employing orthogonal small-scales with the optimality projector. This leads to other well-known stabilization formulations, namely SUPG and GLS.

\subsection{Design condition}
To clarify, let us explicitly mention the design condition of the stabilized formulation which emerges from (\ref{eq:VMS_energy}). We seek for a stabilized formulation corresponding to (\ref{sec:GE, subsec:SF, CD Strong}) which displays local energy behavior as:
\begin{align}\label{sec:EESSMSSS, subsec:EE, energy evo design cond}
  \dfrac{{\rm d}}{{\rm d}t} E_\omega = &-  \| \kappa^{1/2} \nabla \phi^h \|^2_\omega + (\phi' ,f)_{\tilde{\omega}}- (1,F_\omega^h)_{\chi_\omega}\nonumber \\
  &-\|\tau^{-1/2} \phi'\|^2_\omega + (\phi^h,f)_{\omega} .
\end{align}
In this paper we call this \textit{correct-energy behavior}. The positive scalar $\tau$ represents the stabilization parameter of the small-scale equation and equals $\tau=\tau_{\text{static}}$ when using static small-scales as in (\ref{eq:small_stat}).

\subsection{The variational multiscale method with dynamic small-scales}  \label{sec:vms_dyn}

An alternative for replacing the small-scale equation with an algebraic relation is to retain the time-integration and only model the spatial part of the operator. This leads to so-called \textit{dynamic small-scales}, as introduced in \cite{CodPriGuaBad07}. The model equation
\begin{align}\label{eq:tau_dyn}
 \partial_t \hat{\phi}'+\tau_{\text{dyn}}^{-1}\hat{\phi}'= -\mathscr{R}\phi^h,
\end{align}
is now an ordinary differential equation in time. Again, we ignore the $~\tilde{}$ sign in the following. The time derivative in (\ref{eq:tau_dyn}) eliminates the first unwanted temporal part in the energy evolution (\ref{eq:VMS_energy}).
Naturally, the stabilization parameter\footnote{This explains our notation $\tau_{\text{static}}$ in Section \ref{sec:stab} where \textit{static} represents \textit{static small-scales}.} now omits a temporal part, since it is explicitly handled, therefore:
\begin{align}
  \tau_{\text{dyn}}=\left(\tau_{\text{conv}}^{-2}+\tau_{\text{diff}}^{-2}\right)^{-1/2}.
\end{align}
Clearly, the static small-scale equation (\ref{eq:small_stat, der}) does not apply anymore. Therefore, the term $\pd_t \phi'$ is active in the large-scale equation. The VMS stabilized formulation with dynamic small-scales now reads:\\

\textit{Find $\phi^h \in \mathcal{W}^h$ such that for all $w^h \in \WW^h$,}
\begin{subequations}
\label{eq:VMS_weak_dyn}
  \begin{alignat}{3}
  (w^h,\pd_t \phi^h + \pd_t \phi' )_{\Omega}  + (w^h,\B{a} \cdot\nabla \phi^h)_{\Omega} + (\nabla w^h, \kappa\nabla \phi^h)_{\Omega}
   - (\B{a} \cdot\nabla w^h + \kappa \Delta w^h,  \phi')_{\tilde{\Omega}}
  &=& (w^h, f)_{\Omega},  \label{eq:VMS_weak_dyn large} \\
  \pd_t \phi'+\tau_{\text{dyn}}^{-1} \phi'&=& -\mathscr{R}\phi^h. \label{eq:VMS_weak_dyn small}
\end{alignat}
\end{subequations}

To arrive at an energy relation we adopt the same procedure as before. The total local energy linked to this variational form evolves as:
\begin{align}\label{eq:VMS_dyn_energy}
  \dfrac{{\rm d}}{{\rm d}t} E_{\omega}=&-  \| \kappa^{1/2}\nabla \phi^h \|^2_{\omega} + (\phi^h,f)_{\omega}
      - (1,F_\omega^h)_{\chi_\omega} \nonumber  \\
  & -\| \tau_{\text{dyn}}^{-1/2} \phi'\|^2_{\tilde{\omega}}+ (\phi' ,f)_{\tilde{\omega}}+2 (\kappa\Delta \phi^h,\phi')_{\tilde{\omega}},
\end{align}
with $F_\omega^h$ defined in (\ref{energy flux VMS}). We observe that adopting dynamic small-scales indeed eliminates the undesired temporal terms.

\subsection{Orthogonality between the large-scales and the small-scales}
The other unwanted term vanishes when the large-scales and small-scales are appropriately orthogonal with respect to each other, namely
\begin{align}\label{eq:orhto}
(\kappa\Delta \phi^h,\phi')_{\Omega} = 0.
\end{align}
This defines the optimality projector (\ref{eq:opt_proj}) which links the stabilized formulation with the desired energy behavior. Therefore we employ the natural choice for the optimality projector:\\

$\mathscr{P}^h: \phi \in \WW \rightarrow \phi^h \in \WW^h$: \textit{Find $\phi^h \in \mathcal{W}^h$ such that for all $w^h \in \WW^h$,}
\begin{align}
(\kappa\Delta w^h,\phi^h)_{\Omega} = (\kappa\Delta w^h,\phi)_{\Omega}.
\end{align}
This yields the required orthogonality.

\subsection{Consistent SUPG with dynamic small-scales}
Employing the orthogonality (\ref{eq:orhto}) directly in the large-scale equation, leads to the dynamic small-scale version of the
well-known SUPG formulation:\\

\textit{Find $\phi^h \in \mathcal{W}^h$ such that for all $w^h \in \WW^h$,}
\begin{subequations}
\label{sec:CSUPG}
  \begin{alignat}{1}
  (w^h,\pd_t \phi^h+ \pd_t \phi' )_{\Omega}  + (w^h,\B{a} \cdot\nabla \phi^h)_{\Omega} + (\nabla w^h, \kappa\nabla \phi^h)_{\Omega}
   - (\B{a} \cdot\nabla w^h,  \phi')_{\tilde{\Omega}}
  =& (w^h, f)_{\Omega}   \\
  \pd_t \phi'+\tau_{\text{dyn}}^{-1} \phi'=& -\mathscr{R}\phi^h.\label{consistent SUPG small scales}
\end{alignat}
\end{subequations}
Unfortunately, this removes only the contribution from the large-scale equation (\ref{eq:VMS_weak_dyn large}); the contribution of the undesirable term from the small-scale equation (\ref{eq:VMS_weak_dyn small}) remains:
\begin{align}\label{sec:CSUPG_energy}
  \dfrac{{\rm d}}{{\rm d}t} E_{\omega}=&-  \| \kappa^{1/2}\nabla \phi^h \|^2_{\omega} + (\phi^h,f)_{\omega}
      - (1,F_\omega^h)_{\chi_\omega}  \nonumber  \\
  & -\| \tau_{\text{dyn}}^{-1/2} \phi'\|^2_{\tilde{\omega}}+ (\phi' ,f)_{\tilde{\omega}}+ (\kappa\Delta \phi^h,\phi')_{\tilde{\omega}}.
\end{align}
The undetermined sign of the last term indicates that the formulation can still create artificial energy locally.

\subsection{Inconsistent SUPG with dynamic small-scales}

Instead of using the orthogonality (\ref{eq:orhto}) only in the large-scale equation, one can adopt it as well on the small-scales (\ref{eq:small}),(\ref{consistent SUPG small scales}). The resulting SUPG-like  formulation with dynamic small-scales reads:\\

\textit{Find $\phi^h \in \mathcal{W}^h$ such that for all $w^h \in \WW^h$,}
\begin{subequations}
\label{eq:ISUPG}
  \begin{alignat}{1}
  (w^h,\pd_t \phi^h + \pd_t \phi')_{\Omega}  + (w^h,\B{a} \cdot\nabla \phi^h)_{\Omega} + (\nabla w^h, \kappa\nabla \phi^h)_{\Omega}
   - (\B{a} \cdot\nabla w^h,  \phi')_{\tilde{\Omega}}
  =& (w^h, f)_{\Omega}   \\
  \pd_t \phi'+\tau_{\text{dyn}}^{-1}\phi'=& - \pd_t \phi^h - \B{a} \cdot\nabla \phi^h + f.
\end{alignat}
\end{subequations}
The energy evolution linked to this formulation,
\begin{align}\label{eq:ISUPG_energy}
  \dfrac{{\rm d}}{{\rm d}t} E_{\omega}=&-  \| \kappa^{1/2}\nabla \phi^h \|^2_{\omega} + (\phi^h ,f)_{\omega}
      - (1,F_\omega^h)_{\chi_\omega}
  +( \phi' ,f)_{\tilde{\omega}}-\| \tau_{\text{dyn}}^{-1/2} \phi'\|^2_{\tilde{\omega}},
\end{align}
has the desired form. However, this formulation is inconsistent because the small-scales are not forced by a full residual: the diffusive term is removed from the residual. The inconsistent character of the formulation limits the adequacy of this formulation.

\subsection{GLS  with dynamic small-scales (GLSD)}

Another alternative is to use the orthogonality only on the large-scale equation, now with double the magnitude. The diffusive stabilized term does not vanish but flips sign instead. In other words the VMS formulation does not convert to a SUPG formulation but to a GLS one. Hence, the VMS approach with the diffusive optimality projection (\ref{eq:orhto}) leads to the so-called GLSD-statement (the \textit{D} stands for \textit{dynamic}) which reads\\

\textit{Find $\phi^h \in \mathcal{W}^h$ such that for all $w^h \in \WW^h$,}
\begin{subequations}
\label{sec:GLS}
  \begin{alignat}{1}
  (w^h,\pd_t \phi^h+ \pd_t \phi' )_{\Omega}  + (w^h,\B{a} \cdot\nabla \phi^h)_{\Omega} + (\nabla w^h, \kappa\nabla \phi^h)_{\Omega}
   - (\B{a} \cdot\nabla w^h- \kappa \Delta w^h,  \phi')_{\Omega}
  =& (w^h, f)_{\Omega}   \\
 \pd_t \phi'+\tau_{\text{dyn}}^{-1} \phi '=& -\mathscr{R}\phi^h.
\end{alignat}
\end{subequations}
This formulation possesses the desired energy evolution:
\begin{align}\label{sec:GLS_energy}
  \dfrac{{\rm d}}{{\rm d}t} E_{\omega}=&-  \| \kappa^{1/2}\nabla \phi^h \|^2_{\omega} + (\phi^h,f)_{\omega}
      - (1,F_\omega^h)_{\chi_\omega}\nonumber \\
     &-\| \tau_{\text{dyn}}^{-1/2} \phi'\|^2_{\tilde{\omega}}+ (\phi' ,f)_{\tilde{\omega}}.
\end{align}
Comparing with the inconsistent SUPG formulation (\ref{eq:ISUPG}), both variational forms possess the correct-energy behavior. However, this formulation distinguishes itself by its consistent character, i.e. the forcing term in the small-scale equation is driven by the full residual.
\subsubsection*{Remark}
An important observation is: \textit{the GLS formulation is justified in a VMS context by invoking the orthogonality demanded for correct-energy behavior.}

\section{Back to a variational multiscale formulation: a stabilized form with correct-energy evolution}\label{sec:correct2}

Section \ref{sec:correct} justifies with orthogonality arguments a GLS-based formulation depicting correct-energy behavior. That methodology \textit{assumes} orthogonality between the large-scales and small-scales but does not actually \textit{enforce} it. This Section devises an alternative VMS stabilization approach that explicitly \textit{enforces} the required orthogonality.

\subsection{The small-scale solution space}

The weak statements of Section \ref{sec:correct} do not explicitly mention the solution space of the small-scales. The small-scales are effectively pointwise values, i.e. $\phi' : \Omega \times \mathcal{I} \rightarrow \mathbb{R}$. 
Section \ref{sec:stab} reveals that the small-scales live in a properly defined space, that is $\phi'\in \WW'$. The orthogonality projector (\ref{eq:orhto}) leads to the following definition of the small-scale space:
\begin{align}\label{eq:small_space_exp}
  \WW'=\WW'_{H_0^1}:=\left\{\phi \in \WW; \left(\kappa\Delta \eta^h,\phi\right)_{\Omega}=0 \quad \text{for all } \eta^h \in \WW^h\right\}.
\end{align}
Note that the projector is induced by the $H_0^1$-seminorm. This function space enjoys good properties, as indicated in \cite{HugSan06}. The discretization dependence of the stabilization parameter $\tau$ originates from the corresponding restricted Green's function. Consult that paper for details.

Employing the restricted solution space $\WW'$ the dynamic VMS formulation (\ref{eq:VMS_weak_dyn}) subtly modifies to\\

\textit{Find $\phi^h \in \WW^h, \phi'\in \WW'$ such that for all $w^h \in \WW^h$,}
\begin{subequations}
\label{sec:DYN_VMS_RES}
  \begin{alignat}{1}
(w^h,\pd_t \phi^h)_{\Omega}
+ (\nabla w^h, \kappa\nabla \phi^h )_{\Omega}
+ (w^h, \B{a} \cdot\nabla \phi^h )_{\Omega}&\nonumber\\
+(w^h, \pd_t \phi' )_{\tilde{\Omega}}-  (\B{a} \cdot\nabla w^h  + \kappa \Delta w^h ,\phi')_{\tilde{\Omega}}
-(w^h,f)_{\Omega}&  =0,  \\
\pd_t \phi'+\tau_{\text{dyn}}^{-1}\phi'+\mathscr{R}\phi^h & = 0.
\end{alignat}
\end{subequations}
The small-scale solution possesses the correct orthogonality by construction which implies the correct-energy behavior (\ref{sec:GLS_energy}).

However, the restriction of the small-scale solution in the weak form (\ref{sec:DYN_VMS_RES}) is troublesome to directly convert the weak statement into a working numerical method. This is mainly due to the infinite dimensionality of the small-scale space $\WW'$.

\subsection{Enforced orthogonality with a Lagrange multiplier (DO formulation)}

In order to avoid dealing with the restricted subspace (\ref{eq:small_space_exp}), we adopt a Lagrange multiplier setting to reformulate the problem into a mixed formulation. This opens up the search space for $\phi'$, while an explicit constraint is added to satisfy the required orthogonality. A Lagrange multiplier enforces this constraint.
This formulation reads as follows:\\

\textit{Find $\left(\phi^h,\sigma^h\right) \in \WW^h\times\WW^h, \phi': \Omega \times \mathcal{I} \rightarrow \mathbb{R}$ such that for all $(w^h,\eta^h) \in \WW^h\times\WW^h$,}
\begin{subequations}
\label{sec:DYN_VMS_LM}
  \begin{alignat}{2}
  (w^h,\pd_t \phi^h)_{\Omega}+ (\nabla w^h, \kappa\nabla \phi^h )_{\Omega}
  + (w^h, \B{a} \cdot\nabla \phi^h )_{\Omega} &&\nonumber\\
  +(w^h, \pd_t \phi' )_{\tilde{\Omega}}-  (\B{a} \cdot\nabla w^h  + \kappa \Delta w^h ,\phi')_{\tilde{\Omega}} -(w^h,f)_{\Omega}&  =&0,  \\
  \pd_t \phi'+\tau_{\text{dyn}}^{-1}\phi'-\kappa \Delta \sigma^h+\mathscr{R}\phi^h & =& 0,  \\
  \left(\kappa \Delta \eta^h,\phi'\right)_{\tilde{\Omega}}&=&0. \label{sec:DYN_VMS_LM3}
  \end{alignat}
\end{subequations}
We refer to it as \textit{DO} where the \textit{D} and \textit{O} stand for \textit{dynamic} and \textit{orthogonal}, respectively. Here denotes $\sigma^h$ the Lagrange multiplier and $\eta^h$ its associated weighting function.

Note that this formulation asks for $C^1$-continuous basis functions because of the use of second derivatives. This additionally motivates our choice of employing IGA.

\subsection{Local energy evolution of the formulation with enforced orthogonality}
We obtain the energy evolution of the proposed method in a similar fashion as before. Hence, select $w^h=\phi^h$ in the large-scale equation, $\eta^h=\sigma^h + \phi^h$ in the third
equation and multiply the small-scale equation by $\phi'$. Summation of the three equations and reordering leads to:
\begin{align}\label{sec:ALT_energy}
  \dfrac{{\rm d}}{{\rm d}t} E_{\omega}=&-  \| \kappa^{1/2}\nabla \phi^h \|^2_{\omega} + (\phi^h ,f)_{\omega}
      - (1,F_\omega^h)_{\chi_\omega}\nonumber\\
  &-\| \tau_{\text{dyn}}^{-1/2} \phi'\|^2_{\tilde{\omega}}+ (\phi' ,f)_{\tilde{\omega}}.
\end{align}
Note that the terms originating from (\ref{sec:DYN_VMS_LM3}) exactly cancel the undesired orthogonality terms and the small-scale Lagrange multiplier term.

\subsection*{Remark}

The separate energy evolution of the large-scales and small-scales deduces in a similar fashion as above. The energies $E^h_\omega=\tfrac{1}{2}(\phi^h,\phi^h)_{\omega}$ and $E'_\omega=\tfrac{1}{2}(\phi',\phi')_{\tilde{\omega}}$ evolve respectively as
\begin{subequations}
\label{eq:energy backscatter}
\begin{alignat}{1}
  \dfrac{{\rm d}}{{\rm d}t} E^h_\omega & = - \| \kappa^{1/2} \nabla \phi^h \|^2_\omega + (\phi^h,f)_\omega+  (\B{a} \cdot \nabla \phi^h ,\phi' )_{\tilde{\omega}}  - (\phi^h,\pd_t \phi' )_{\tilde{\omega}} - (1,F_\omega^h)_{\chi_\omega}, \label{eq:energy backscatter: large}\\
  \dfrac{{\rm d}}{{\rm d}t} E'_\omega  & = - \| \tau_{\text{dyn}}^{-1/2} \phi'\|^2_{\tilde{\omega}}   + (\phi',f)_{\tilde{\omega}}-  ( \B{a} \cdot \nabla \phi^h ,  \phi')_{\tilde{\omega}}  -  ( \phi', \pd_t \phi^h )_{\tilde{\omega}}. \label{eq:energy backscatter: small}
\end{alignat}
\end{subequations}
The first term of (\ref{eq:energy backscatter: small}) may be viewed as diffusion of the small-scales. The convective contributions exchange energy between the large-scales and small-scales. It is important to emphasize that these energies do not add up to the total local energy $E_\omega$: the cross terms are missing. Their contributions appear in both (\ref{eq:energy backscatter: large})-({\ref{eq:energy backscatter: small}).

\section{Temporal-integration}\label{sec:time}

This Section is devoted to the time-integration for which we adopt the generalized-$\alpha$ time integrator. We start off with a brief recap of the method in a general setting, after which we discuss the use of this method for the small-scales particularly. The reminder presents the collection of time-integrators within this framework which yields a concrete energy evolution statement of consecutive time levels.

\subsection{The generalized-$\alpha$ time integrator}\label{subsec:gen_alpha}

We employ the generalized-$\alpha$ method for the temporal-integration \cite{Hul93}. The algorithm
reads: \\

\textit{Given the data $\phi_n,\dot{\phi}_n$, find $\phi_{n+\alpha_f}, \dot{\phi}_{n+\alpha_m}, \phi_{n+1}, \dot{\phi}_{n+1}$ such that }
\begin{subequations}\label{generalized alpha method}
\label{eq:gen_alpha}
\begin{alignat}{1}
  \dot{\phi}_{n+\alpha_m} = & \mathcal{G}(\phi_{n+\alpha_f}) \label{ODE gen alpha},\\
  \phi_{n+1}               = & \phi_{n}  + \Delta t \left((1 - \gamma) \dot {\phi}_{n}  +  \gamma \dot {\phi}_{n+1}\right)\label{sec:EECDE, subsec:GESWF, TI, phi 1},\\
  \dot {\phi}_{n+\alpha_m} = & (1 - \alpha_m)\dot {\phi}_{n} + \alpha_m\dot {\phi}_{n+1} \label{sec:EECDE, subsec:GESWF, TI, phi alpham},\\
  \phi_{n+\alpha_f}        = & (1 - \alpha_f)\phi_{n} + \alpha_f\phi_{n+1}.\label{sec:EECDE, subsec:GESWF, TI, phi alphaf}
\end{alignat}
\end{subequations}
Here $\pd_t \phi=\mathcal{G}(\phi)$ is the governing ordinary differential equation, $\dot{\phi}$ is the
discretized time derivative and the time step size is $\Delta t=t_{n+1}-t_n$.
The scalars $\alpha_f, \alpha_m, \gamma$ are algorithmic coefficients that need to be properly selected.
 The methods reduce to some of the classical time integrators for specific choices of the
  time-step parameters, e.g. for $\alpha_f=\alpha_m=\gamma=1$ to backward Euler
  and for $\alpha_f=\alpha_m=\gamma=\frac{1}{2}$ to Crank-Nicolson.
   It is unconditionally stable if $\alpha_m\geq \alpha_f \geq \frac{1}{2}$ (i.e. when it is more implicit than explicit).
   The second-order accuracy requirement reads \cite{Hul93,WhiJan99}:
\begin{align}
  \label{sec:EECDE, subsec:GESWF, TI, 2nd acc}
  \gamma = \frac{1}{2} +\alpha_m -\alpha_f.
\end{align}

\subsection{Time-integration of the small-scales}
\label{sec:time, subsec:TIFSE}
The small-scale equations are ordinary differential equations. Employing (\ref{eq:tau_dyn}) for (\ref{ODE gen alpha}) an explicit solution of system (\ref{eq:gen_alpha}) directly follows
\begin{align}\label{eq:ga_small}
    \dot{\phi}'_{n+1}=&C\left(-\frac{1}{\gamma\Delta t}\dot{\phi}'_{n}\left(1-\alpha_m+(1-\gamma)\alpha_f \frac{\Delta t}{\tau_{\text{dyn}}} \right)-\frac{1}{\tau_{\text{dyn}}\gamma\Delta t}\phi'_n-\frac{\mathscr{R}_{n+\alpha}^h}{\gamma\Delta t}\right),
\end{align}
with constant $C=\alpha_f^{-1}\left(\tau_{\text{time}}^{-1}+\tau_{\text{dyn}}^{-1}\right)^{-1}$.

When using dynamic small-scales, the stabilizing properties of the weak formulation depend on the relation between the small- and the large-scales. This relation also enters in the Jacobian of the weak formulation. To this purpose we now explore this link. Let us define the so-called \textit{effective stabilization parameter} as follows
\begin{align}
\tau_{\text{eff}} =& - \frac{\pd \phi'_{n+\alpha_f} }{\pd \mathscr{R}_{n+\alpha}^h},
\end{align}
inspired by (\ref{eq:small_stat, algebraic}). In the case of static small-scales, depicted in (\ref{eq:small_stat}), the trivial expression yields
\begin{align}
\tau_{\text{eff}}  = \tau_{\text{stat}}=\left(\tau_{\text{time}}^{-2}+ \tau_{\text{dyn}}^{-2}\right)^{-1/2}. 
\end{align}
When employing dynamic subscales as in (\ref{eq:ga_small}), we get
\begin{align}
\tau_{\text{eff}} =& - \frac{\pd \phi'_{n+\alpha_f} }{\pd \mathscr{R}_{n+\alpha}^h}
       = - \frac{\pd \phi'_{n+\alpha_f} }{\pd \phi'_{n+1}  }
         \cdot\frac{\pd \phi'_{n+1} }{\pd \dot{\phi}'_{n+1} }
          \cdot \frac{\pd \dot{\phi}'_{n+1} } {\pd \mathscr{R}_{n+\alpha}^h} \nonumber \\
        =&
\alpha_f
 \cdot\gamma\Delta t
 \cdot  \frac{C}{\gamma\Delta t}
        =
 \left(\alpha_m \alpha_f^{-1} \gamma^{-1} \Delta t^{-1}+ \tau_{\text{dyn}}^{-1}\right)^{-1}
 = \left(\tau_{\text{time}}^{-1}+ \tau_{\text{dyn}}^{-1}\right)^{-1},
\end{align}
from which our definition of $\tau_{\text{time}}$ is inspired:
\begin{equation}
 \tau_{\text{time}} := \dfrac{\alpha_f \gamma \Delta t}{\alpha_m}.
\end{equation}
The effective stabilization parameter $\tau_{\text{eff}}$ is very similar to $\tau_{\text{stat}}$ and shows the same asymptotic behavior. This modification of stabilization parameter effectively constitutes a change in the so-called \textit{r-switch} \cite{TeOs00} from $r=2$ to $r=1$. The \textit{r-switch} is a smooth approximation of the minimum operator. A high value of the integer $r$ indicates a sharp switch. In fact, when the stabilization parameters are defined with the r-switch of $r=1$:
\begin{subequations}
\label{eq:tau_mod}
\begin{alignat}{1}
  \widetilde \tau_{\text{stat}}=&\left(\tau^{-1}_{\text{conv}}+\tau^{-1}_{\text{diff}} + \tau^{-1}_{\text{time}}\right)^{-1},\\
  \widetilde \tau_{\text{dyn}}=&\left(\tau^{-1}_{\text{conv}}+\tau^{-1}_{\text{diff}}\right)^{-1},
\end{alignat}
\end{subequations}
the effective stabilization parameters would be identical.

\subsection{Proper energy evolution}\label{sec:TI, subsec:TILED}

The energy evolution equations (\ref{sec:GLS_energy}) or (\ref{sec:ALT_energy}) reveal a (global) guaranteed energy decay in the absence of external forcing and boundaries, namely,
\begin{align}
\dfrac{{\rm d}}{{\rm d}t} E_\Omega =  - \|\kappa^{1/2}\nabla \phi^h \|^2_{\Omega} -  \|\tau_{\text{dyn}}^{-1/2}  \phi' \|^2_{\tilde{\Omega}}.
\end{align}
The time-integrator should obey this decaying property. Moreover, ideally it leads to a guaranteed decay of energy for consecutive time levels, that is,
\begin{align}
 E_{n+1}\leq E_n.
\end{align}
To arrive at an appropriate energy statement when employing the generalized-$\alpha$ method the procedure reads as follows. Multiply the small-scale equation with $\phi'_{n+\alpha_f}$, integrate the result and add it to the weak form in which  $w^h=\phi^h_{n+\alpha_f}$ is selected. The continuous form of this approach has been demonstrated before in this paper, see e.g. Section \ref{sec:stab, subsec: Local energy evolution of the VMS formulation}. This leads to the correct
symmetric inner products for the spatial terms, and proper norms therefore.
Additionally, the temporal terms, leading to the energy derivatives yield,
\begin{align}\label{energy time derivatives}
\Delta E = \Delta t \dot{E}_{n+\alpha} = \Delta t (\phi_{n+\alpha_f} ,\dot {\phi}_{n+\alpha_m} )_{\Omega},
\end{align}
where $\dot{E}_{n+\alpha}$ is the temporal derivative of the energy at time level $n+\alpha$ and $\phi=\phi^h+\phi'$. The abuse of notation demands the integration to be interpreted on $\tilde{\Omega}$ for terms containing the small-scales. In the following we derive time-stepping parameters within the generalized-$\alpha$ time integrator framework which link this temporal term to a proper energy behavior.

Substitution of (\ref{sec:EECDE, subsec:GESWF, TI, phi 1})-(\ref{sec:EECDE, subsec:GESWF, TI, phi alphaf}) into (\ref{energy time derivatives}) yields:
\begin{align}
\Delta t \dot{E}_{n+\alpha} = \nonumber&
 \Delta t(\phi_{n+\alpha_f} ,\dot {\phi}_{n+\alpha_m} )_{\Omega} \nonumber\\
= &  \Delta t( (1 - \alpha_f)\phi_{n} + \alpha_f\phi_{n+1},
      (1 - \alpha_m)\dot{\phi}_{n} + \alpha_m\dot{\phi}_{n+1} )_{\Omega} \nonumber \\
= &  ( (1 - \alpha_f)\phi_{n} + \alpha_f\phi_{n+1},
   \left (1 -    \frac{\alpha_m}{\gamma}   \right )\Delta t \dot {\phi}_{n} +  \frac{\alpha_m}{\gamma} \left ( \phi_{n+1}  - \phi_{n}  \right ) )_{\Omega} \nonumber \\
= & - \frac{(1 - \alpha_f)\alpha_m}{\gamma }(\phi_{n}   ,\phi_{n} )_{\Omega}
    + \frac{\alpha_f\alpha_m}{\gamma}  (\phi_{n+1}   ,\phi_{n+1} )_{\Omega}
    + (1 - 2\alpha_f ) \frac{\alpha_m}{ \gamma}  (\phi_{n}   ,\phi_{n+1} )_{\Omega}\nonumber \\
 &  + \left ( (1 - \alpha_f)\phi_{n} + \alpha_f\phi_{n+1},
      \left (1 -    \frac{\alpha_m}{\gamma}   \right )\Delta t \dot {\phi}_{n} \right)_{\Omega},
\end{align}
The last term is precarious. The sign of the temporal derivative $\dot {\phi}_{n}$ is not determined. It appears without $(n+1)$-counterpart which leads to an uncontrollable last term. We remedy this issue by requiring the last term to vanish. This occurs when $\alpha_m = \gamma$. Following this road, the temporal term becomes
\begin{align}
  \Delta t \dot{E}_{n+\alpha} =&  \alpha_f (\phi_{n+1}   ,\phi_{n+1} )_{\Omega}
                        - (1 - \alpha_f)(\phi_{n}   ,\phi_{n} )_{\Omega}
                        + (1 - 2\alpha_f ) (\phi_{n}   ,\phi_{n+1} )_{\Omega} \nonumber \\
                        =&  E_{n+1}-E_n+\left(\alpha_f-\tfrac{1}{2}\right)\left[\left(\phi_{n+1},\phi_{n+1}\right)_{\Omega}-2\left(\phi_{n+1},\phi_{n}\right)_{\Omega}+\left(\phi_{n},\phi_{n}\right)_{\Omega}\right] \nonumber \\
                        =&  E_{n+1}-E_{n}+ (\alpha_f-\onehalf) \| \phi_{n+1}  - \phi_{n} \|^2_{\Omega} \nonumber \\
                        =&  E_{n+1}-E_{n}+ \Delta t^2(\alpha_f-\onehalf) \| \dot{\phi}_{n+\alpha_m} \|^2_{\Omega}.
\end{align}
where the last equality is a direct consequence of (\ref{generalized alpha method}) with $\alpha_m = \gamma$. Henceforth, by combining this equation with (\ref{sec:ALT_energy}) the discretized energy (of the DO form) satisfies
\begin{align}\label{sec:discretized energy VMSDO}
  \dfrac{E_{n+1}-E_n}{\Delta t}+\Delta t (\alpha_f-\onehalf) \| \dot{\phi}_{n+\alpha_m} \|^2_{\Omega}=&-  \| \kappa^{1/2}\nabla \phi^h_{n+\alpha_f} \|^2_{\Omega} + (\phi^h_{n+\alpha_f} ,f)_{\Omega}\nonumber \\
  &-\| \tau_{\text{dyn}}^{-1/2} \phi_{n+\alpha_f}'\|^2_{\tilde{\Omega}}+ (\phi_{n+\alpha_f}' ,f)_{\tilde{\Omega}}.
\end{align}
The trivially equivalent form
\begin{align}\label{sec:discretized energy VMSDO equivalent}
  E_{n+1}=E_n-\Delta t^2 (\alpha_f-\onehalf) \| \dot{\phi}_{n+\alpha_m} \|^2_{\Omega}&-  \Delta t\| \kappa^{1/2}\nabla \phi^h_{n+\alpha_f} \|^2_{\Omega}-\Delta t\| \tau_{\text{dyn}}^{-1/2} \phi_{n+\alpha_f}'\|^2_{\tilde{\Omega}} \nonumber \\
  &+ \Delta t(\phi^h_{n+\alpha_f} ,f)_{\Omega}+ \Delta t(\phi_{n+\alpha_f}' ,f)_{\tilde{\Omega}}
\end{align}
reveals that a decay of the discretized energy is guaranteed when, in absence of forcing, additionally $\alpha_f \geq \tfrac{1}{2}$. The first term on the right-hand side, which again should be interpreted on $\tilde{\Omega}$ for the small-scales, is numerical diffusion which vanishes for $\alpha_f=\tfrac{1}{2}$. Hence, the parameter family $\alpha_f\geq\tfrac{1}{2}, \alpha_m=\gamma$, which includes the Crank-Nicolson time-integrator, can be linked to a proper energy decay. Notice that for $\alpha_f=\tfrac{1}{2}$ the stability constraint is fulfilled and the second-order accuracy condition (\ref{sec:EECDE, subsec:GESWF, TI, 2nd acc}) is not harmed.
\clearpage
\section{Numerical verification}
\label{sec:Numerical verification}
This Section provides the numerical verification of the proposed variational formulations of the Sections \ref{sec:stab}-\ref{sec:correct2} on a model problem. We focus on the energy behavior on both a global and a local level. First, the energy behavior is assessed verifying the overall performance of the newly proposed methods. Next, we zoom in on the effect of the small-scales on the energy dissipation. 
We assess its global evolution and local distribution and examine the contributions of the unwanted terms.

\subsection{Model problem description}
The problem under consideration is convection skew to the mesh on a $1x1$-domain with periodic boundaries.
The velocity is $\B{a}=(1,1)$, therefore the profile loops once through the mesh and arrives at its start position at $t=1.0$.
The diffusivity is set to $\kappa = 5\times 10^{-4}$. No external forcing is applied.
The initial condition is a sharp block of the form:
\begin{align}\label{model problem: IC0}
  \phi_0(\bx)  & =H(|x-\tfrac{1}{2}|)H(|y-\tfrac{1}{2}|),\\
  H(z) &= \begin{cases}
	    1                   & \quad z<l\\
            0                   & \quad l<z,
  \end{cases} \label{sec:NV, IC}
\end{align}
where $l$ is a specified length.
For the discretization we employ NURBS\footnote{Note that the quadratic NURBS reduce to B-splines on our uniform Cartesian mesh.}.
Linear NURBS are not considered as they would eliminate the diffusive stabilization term $(\kappa \Delta \phi^h,\phi')_\Omega$ and hence the stabilized forms (SUPG, VMS and GLS) coincide.
All our implementations use quadratic NURBS to bypass this effect.  The sharpest profile that can be exactly represented on the mesh has the form:
\begin{align}\label{model problem: IC}
  \phi_0(\bx)  &= \hat{H}(|x-\tfrac{1}{2}|)\hat{H}(|y-\tfrac{1}{2}|),\\
  \hat{H}(z) &= \begin{cases}
	    1                   & \quad \quad \quad ~ z<l_0\\
            1-\frac{(z-l_0)^2}{2h_c^2} & \quad l_0<z<l_1\\
            \frac{(l_2-z)^2}{2h_c^2}   & \quad l_1<z<l_2\\
            0                   & \quad l_2<z,
  \end{cases} \label{sec:NV, IC}
\end{align}
where $l_0$, $l_1$ and $l_2$ are specified lengths of the different segments that have to coincide with mesh lines. 
We will use $16$x$16$,  $32$x$32$ and  $64$x$64$-element meshes. 
As we want to verify the behavior of the method itself and not consider the error in representing the initial condition we use the exact same initial condition on all meshes.
This is in this case the initial condition of the $16$x$16$-element mesh. 
Therefore we choose  $l_0=nh_c,~l_1=(n+1)h_c$ and $l_2=(n+2)h_c$ with  $n=2$ and $h_c=\tfrac{1}{16}$. 

The implementations use the energy-conserving time-integrator of Section \ref{sec:time} with $\alpha_f=\tfrac{1}{2}$ motivated by both the second-order temporal accuracy and the stability. All computations use a CFL number of $0.5$.
\begin{figure}[h!]
    \begin{center}
    \vspace{1cm} \begin{subfigure}[b]{0.45\textwidth}
        \includegraphics[scale=0.1]{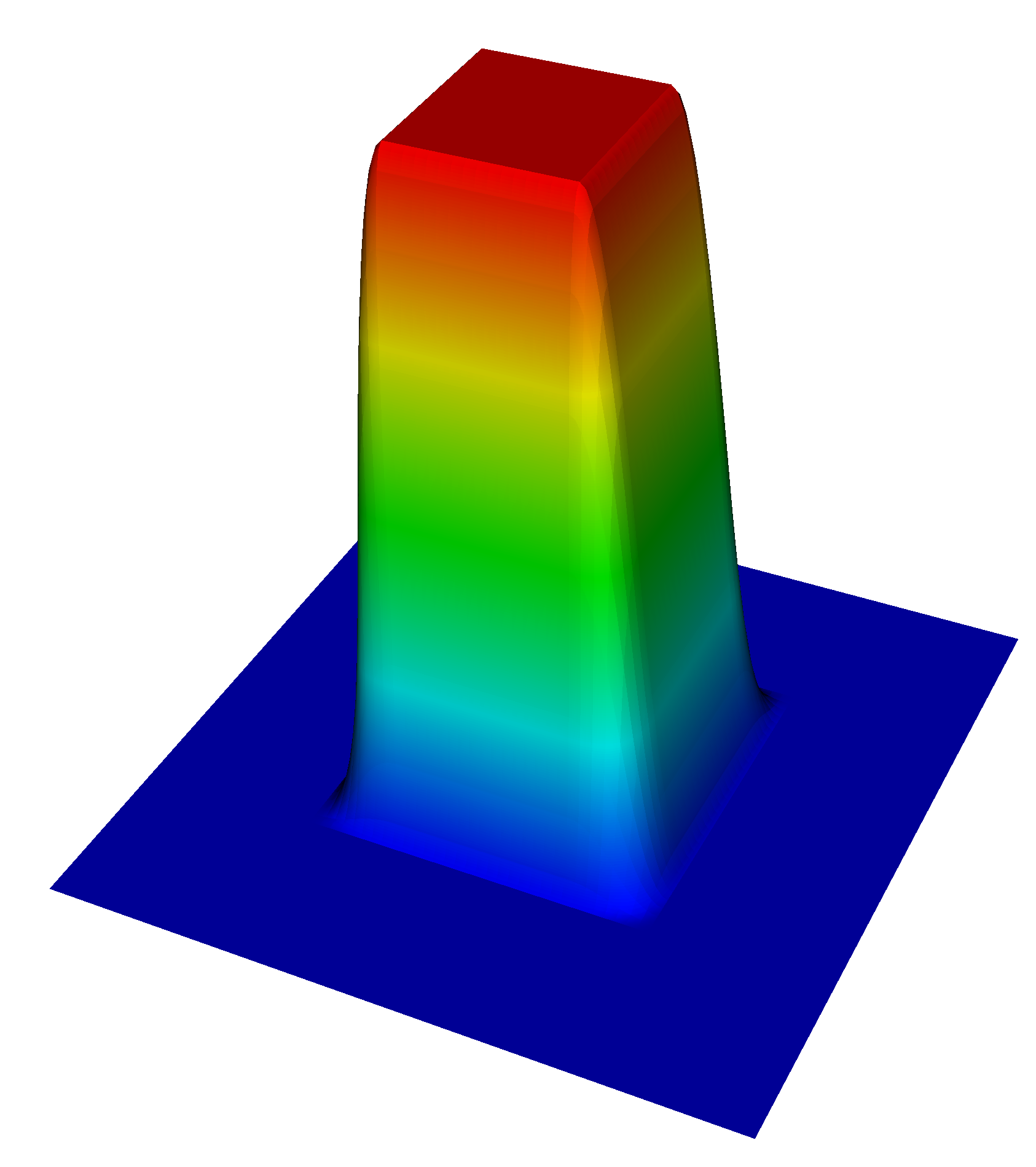}
        \caption{$t=0.0$}
    \end{subfigure}
    \begin{subfigure}[b]{0.45\textwidth}
        \includegraphics[scale=0.1]{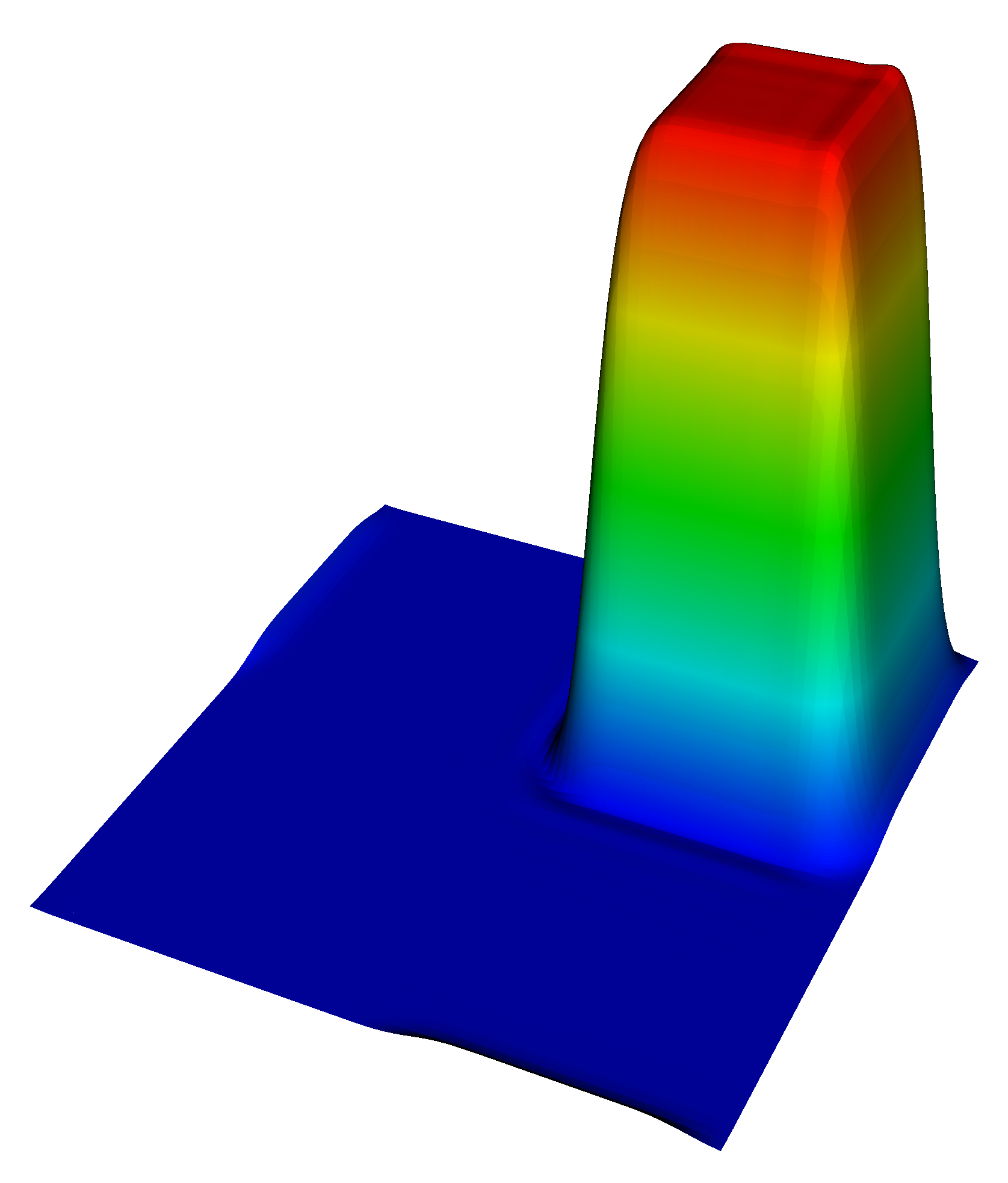}
        \caption{$t=0.25$}
    \end{subfigure}
    \begin{subfigure}[b]{0.45\textwidth}
        \includegraphics[scale=0.1]{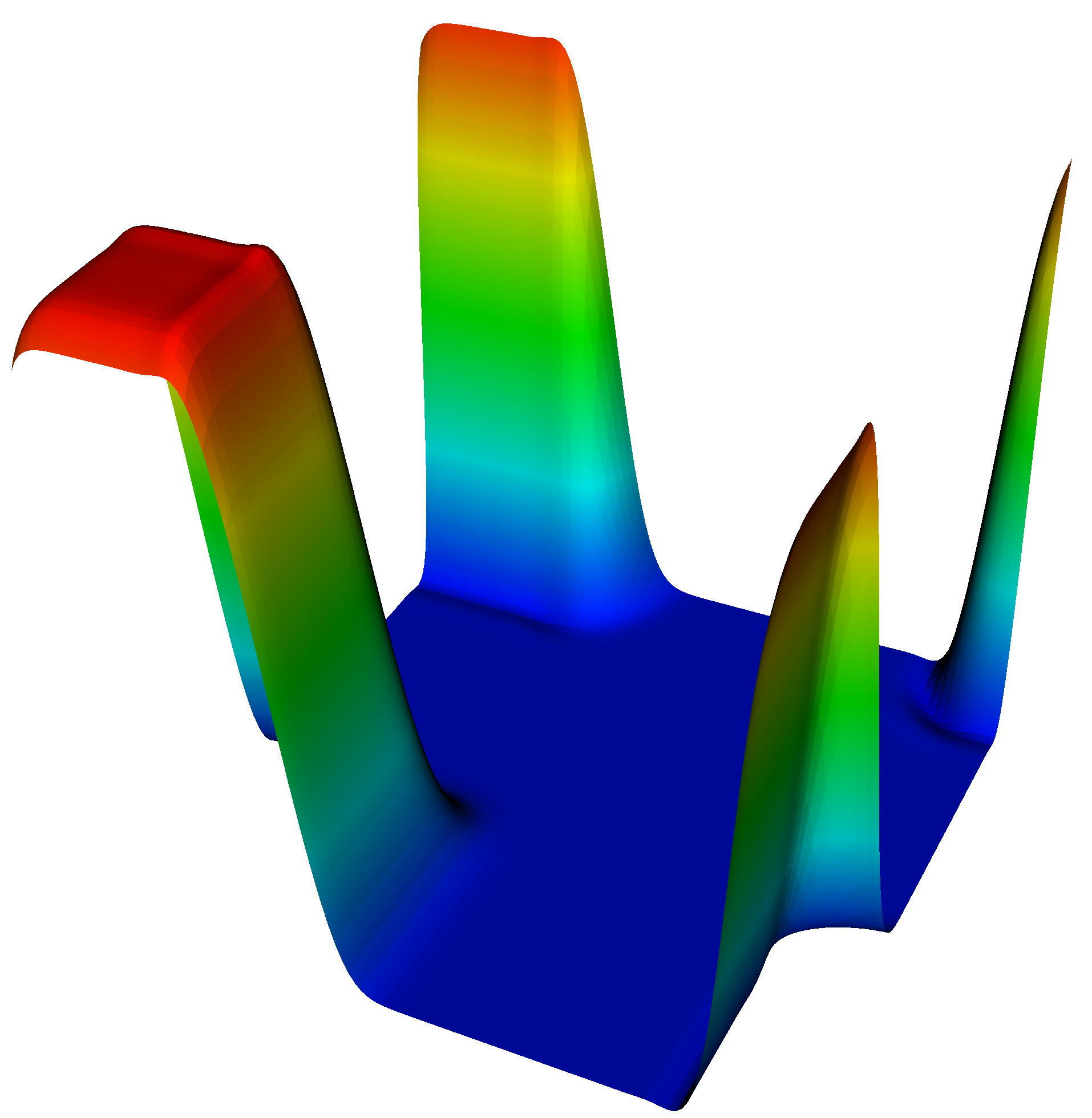}
        \caption{$t=0.625$}
    \end{subfigure}
    \begin{subfigure}[b]{0.45\textwidth}
        \includegraphics[scale=0.1]{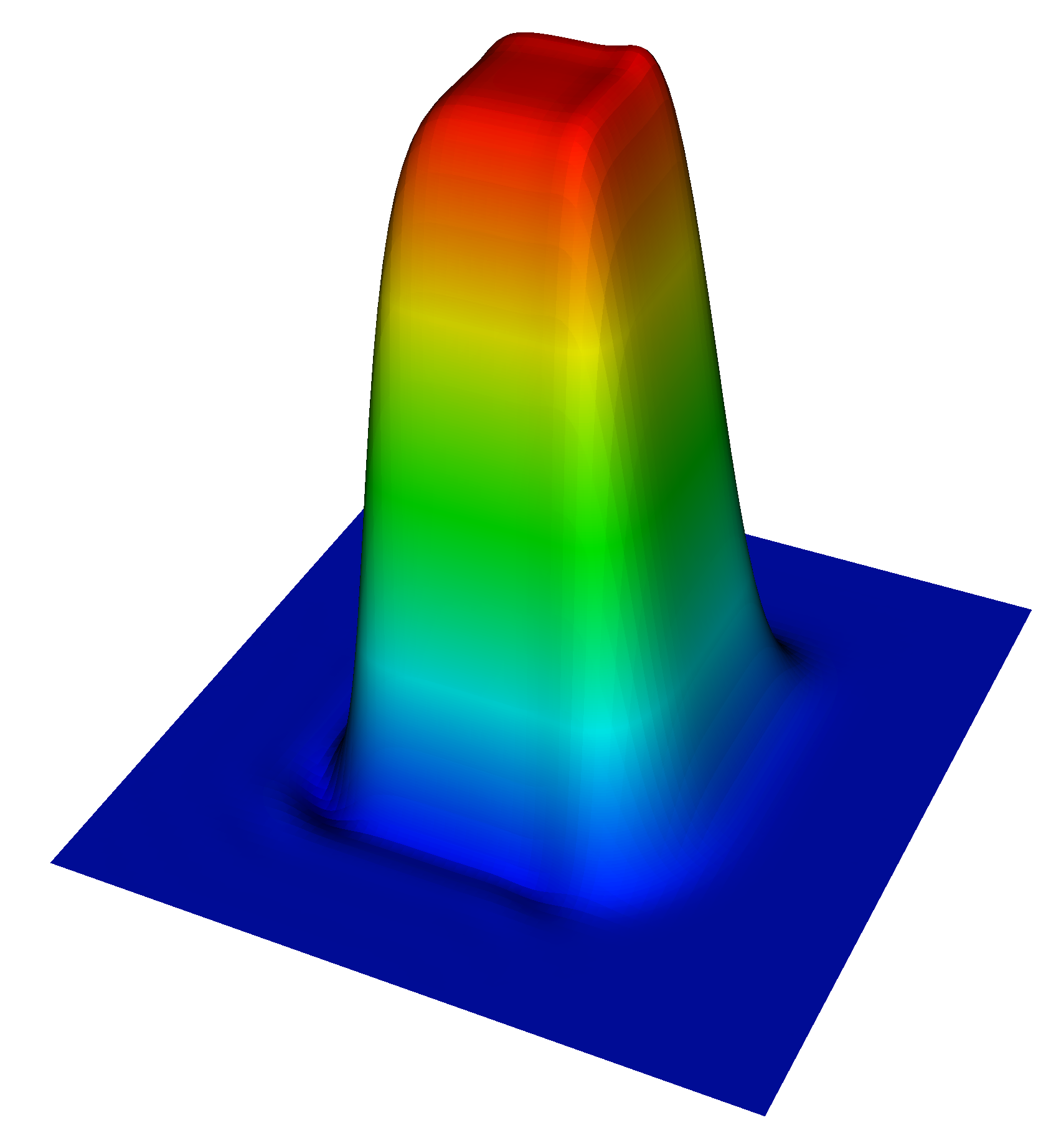}
        \caption{$t=1.0$}
    \end{subfigure}
    \caption{Time evolution of the block profile convected through the mesh.}\label{fig:time evolution}
    \end{center}
\end{figure}
Figure \ref{fig:time evolution} shows the profile traveling through the mesh from $t=0$ until $t=1.0$. The profile exits the mesh approximately halfway during the simulation (at $t=0.5$) and enters at the opposite corner due to the periodic boundary conditions.

In the following we present energy evolution results for three different methods: (i) the SUPG method with static small-scales (SUPGS), (ii) the GLS method with dynamic small-scales (GLSD) and (iii) the dynamic orthogonal formulation (DO). 
These were chosen because the last two exhibit the correct-energy behavior, while SUPG with static small-scales is the classical approach and serves as a reference. It turns out that all methods with static small-scales show very similar behavior.

\clearpage

\subsection{Overall energy behavior}
Figure \ref{conv:SUPGS,GLSD,DO} displays the energy behavior for various methods on different meshes. It shows convergence of the energy evolution for each one of the methods. 
For the SUPGS we have two alternative energy definitions, i.e. one based on only the large-scales and one based on both large- and small-scales (denoted as total energy). 
The energy behavior on the 16x16-mesh is not converged yet whereas the energy on the 32x32-mesh already closely follows that of the finer meshes. In the following we study in more detail the energy evolution on a 32x32-mesh. At this stage there is no visible difference between these solutions.
\begin{figure}[h!]
    \centering
    \begin{subfigure}[b]{0.45\textwidth}
        \includegraphics[scale=0.60]{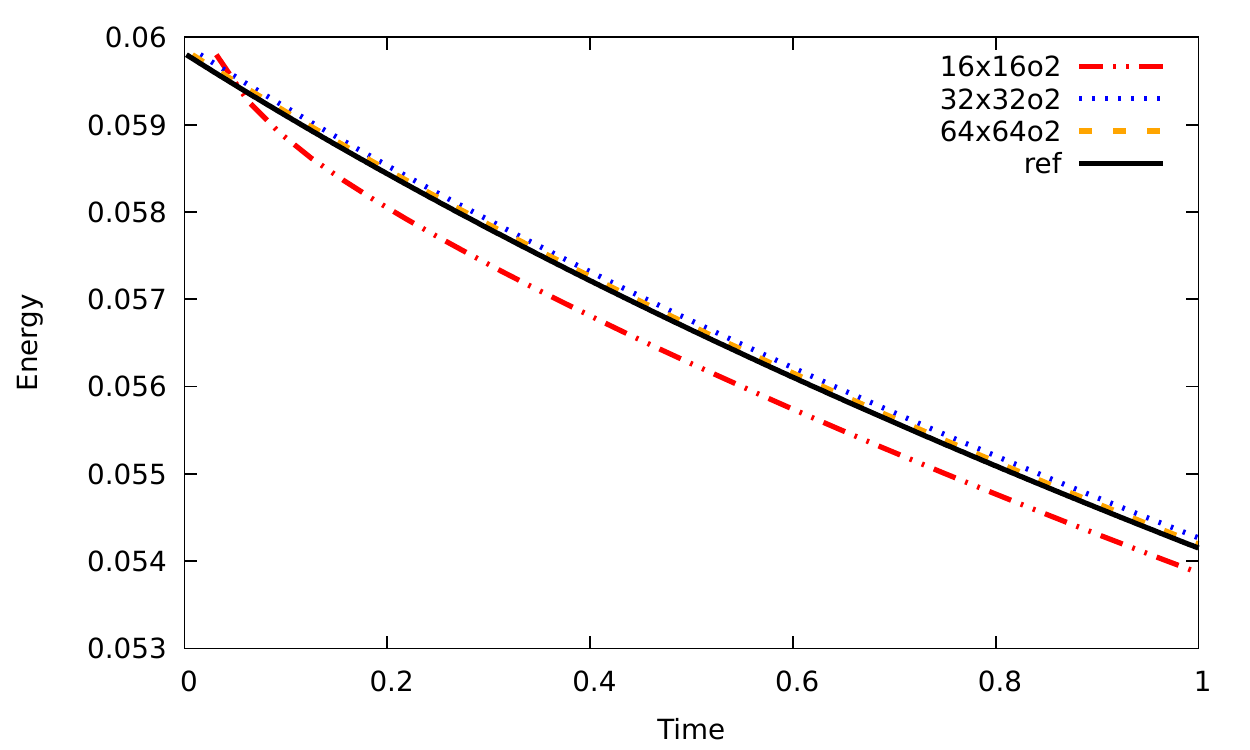}
        \caption{SUPGS, large-scale energy}
    \end{subfigure}
    \quad
    \begin{subfigure}[b]{0.45\textwidth}
        \includegraphics[scale=0.60]{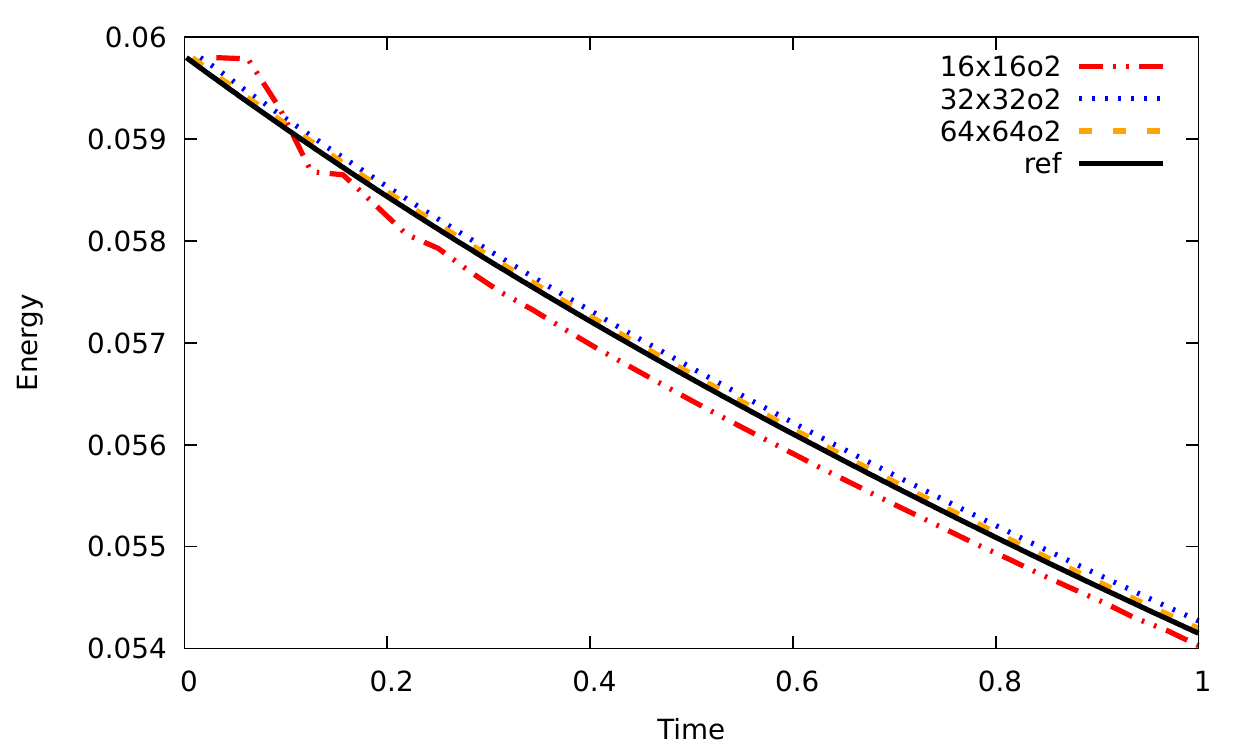}
        \caption{SUPGS, total energy}
    \end{subfigure}
    \begin{subfigure}[b]{0.45\textwidth}
        \includegraphics[scale=0.60]{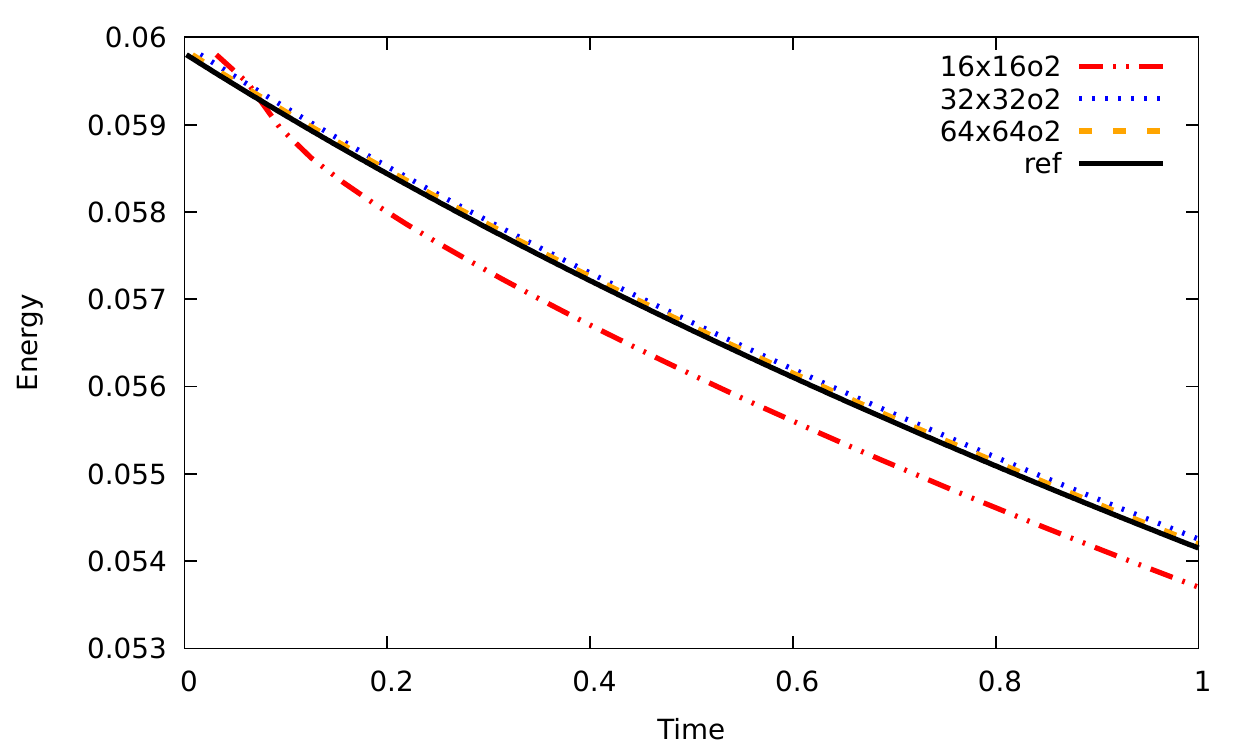}
        \caption{GLSD, total energy}
    \end{subfigure}
    \quad
    \begin{subfigure}[b]{0.45\textwidth}
        \includegraphics[scale=0.60]{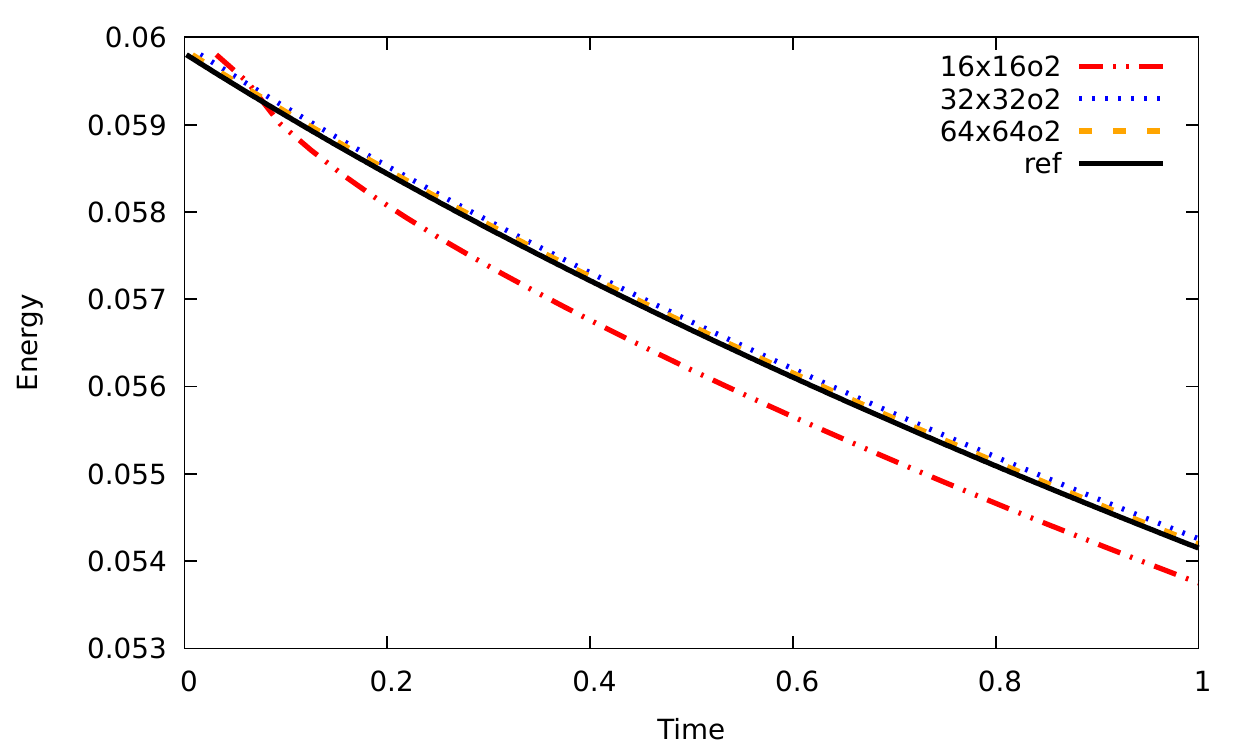}
        \caption{DO, total energy}
    \end{subfigure}
    \caption{Energy evolution for various meshes: (a) the large-scale energy for SUPGS and the total energy for (b) SUPGS, (c) GLSD and (d) DO. An overkill reference solution is added (the continuous black line).}\label{conv:SUPGS,GLSD,DO}
\end{figure}\\

\clearpage
\subsection{Energy dissipation by the small-scales}
Here we study the effect of the small-scales on the energy dissipation. The choice $\alpha_f=\alpha_m=\gamma=\tfrac{1}{2}$ removes the effect of the time-integrator on the energy dissipation. The energy evolution for SUPGS takes the form:
\begin{subequations}
\label{eq:SUPG evolution}
\begin{alignat}{1}
  \dfrac{{\rm d}}{{\rm d}t}E^h_{\omega}=&-\| \kappa^{1/2}\nabla \phi^h \|^2_{\omega}- (1,F_\omega^h)_{\chi_\omega}-\| \tau_{\text{stat}}^{-1/2} \phi'\|^2_{\tilde{\omega}}+(\kappa \Delta \phi^h,\phi')_{\tilde{\omega}}-(\phi',\pd_t \phi^h)_{\tilde{\omega}},\\
  \dfrac{{\rm d}}{{\rm d}t}E_{\omega}=&-\| \kappa^{1/2}\nabla \phi^h \|^2_{\omega}- (1,F_\omega^h)_{\chi_\omega}-\| \tau_{\text{stat}}^{-1/2} \phi'\|^2_{\tilde{\omega}}+(\kappa \Delta \phi^h,\phi')_{\tilde{\omega}}+(\pd_t \phi',\phi^h+\phi')_{\tilde{\omega}}
\end{alignat}
\end{subequations}
for the large-scale energy and the total energy respectively. The GLSD method and the DO formulation show correct-energy evolution:
\begin{align}\label{energy evolution gsld do}
  \dfrac{{\rm d}}{{\rm d}t}E_{\omega}=&-\| \kappa^{1/2}\nabla \phi^h \|^2_{\omega}- (1,F_\omega^h)_{\chi_\omega}-\| \tau_{\text{dyn}}^{-1/2} \phi'\|^2_{\tilde{\omega}}.
\end{align}
The right-hand side terms are evaluated at time level $n+1/2$. The last three terms of each of (\ref{eq:SUPG evolution}) and the last term of (\ref{energy evolution gsld do}) represent the small-scale contribution to the energy dissipation. Figures \ref{fig:Small-scale dissipation global} displays the evolution of the small-scale contribution to energy dissipation on a global scale. 

\begin{figure}[h!]
    \begin{center}
    \begin{subfigure}[b]{0.45\textwidth}
        \includegraphics[scale=0.55]{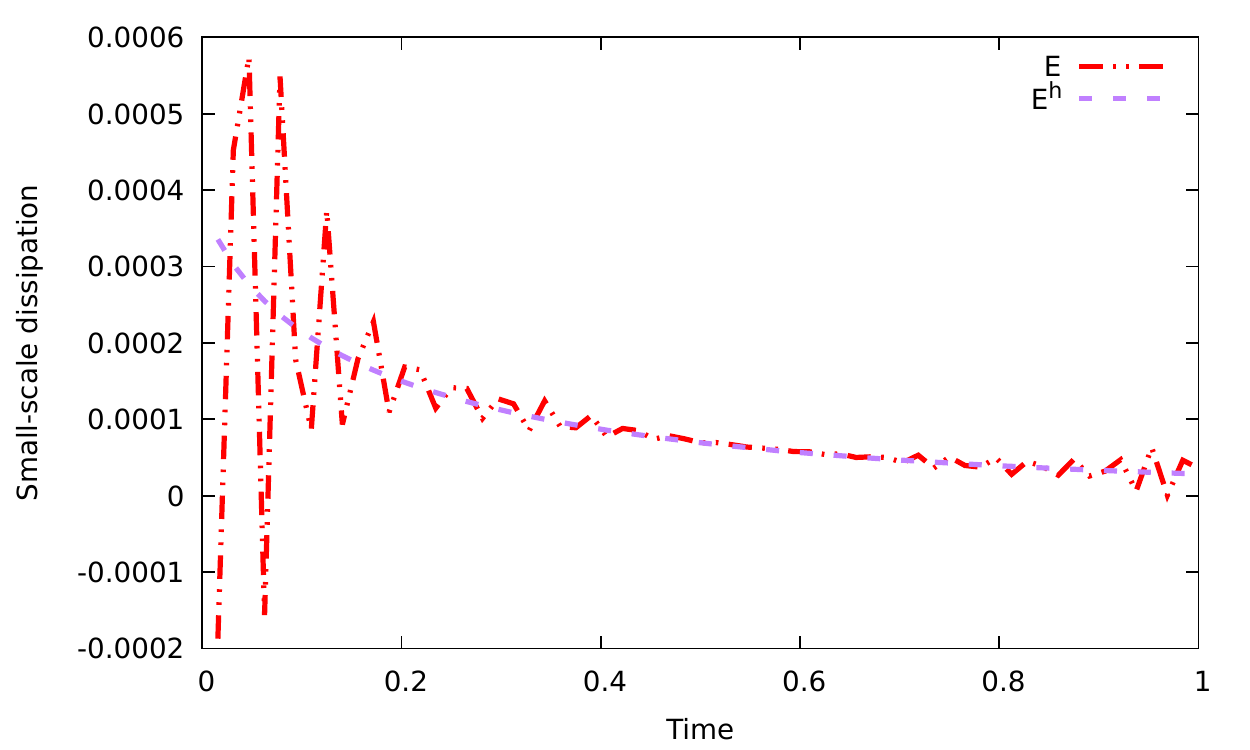}
        \caption{SUPGS}
    \end{subfigure}
    \quad   
    \begin{subfigure}[b]{0.45\textwidth}
        \includegraphics[scale=0.55]{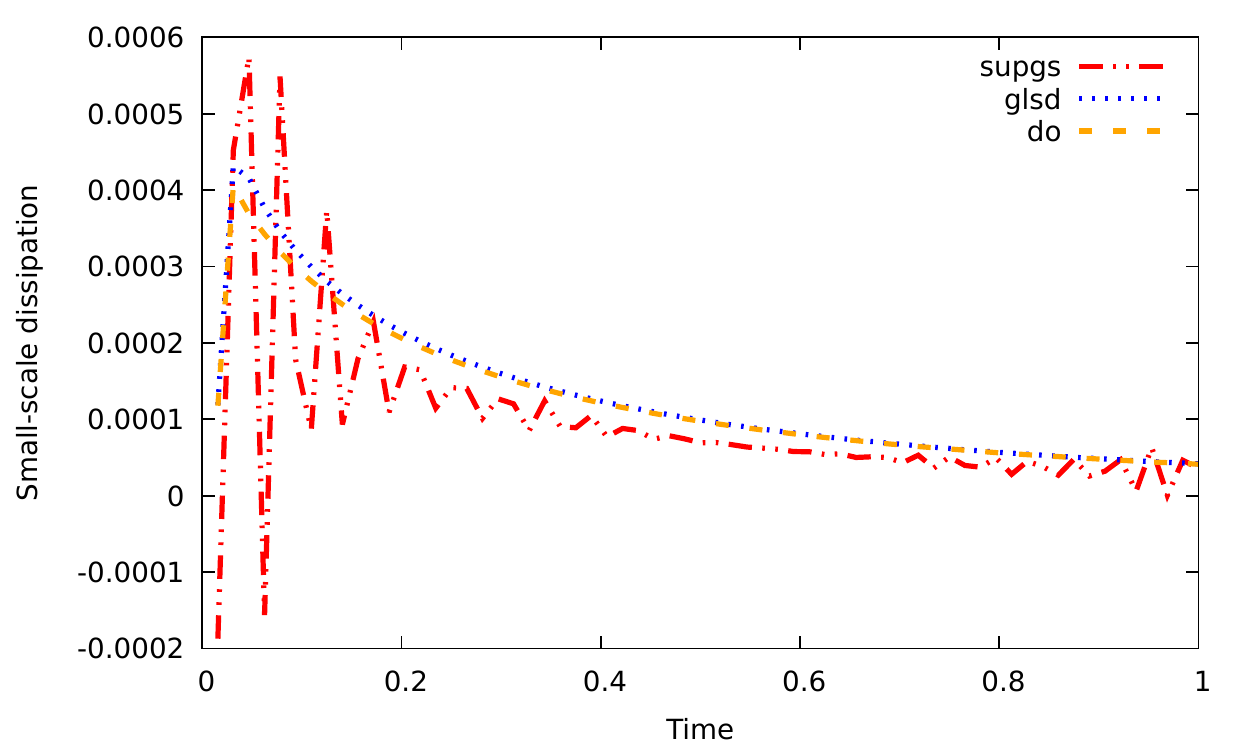}
        \caption{Total energy evolution contribution}
    \end{subfigure}
    \caption{Evolution of the small-scale contribution to energy dissipation on a global scale: (a) large-scale ($E^h$) and total energy ($E$) for SUPG with static small-scales, (b) total energy
for the three methods.}
    \label{fig:Small-scale dissipation global}
    \end{center}
\end{figure}
As anticipated GLSD and DO show positive energy dissipation. On the other hand it is clear that, when considering the total energy, SUPGS has problematic dissipation behavior. It shows severe wiggles resulting in undershoots with negative dissipation. However, when considering large-scale energy there seems to be no problem. 

Figure \ref{fig:Small-scale dissipation local}  shows a typical  local distribution of the small-scale dissipation. These largely confirm the findings from Figure \ref{fig:Small-scale dissipation global}.
GLSD and DO show strictly positive energy dissipation throughout the domain. For SUPGS now both energy definitions show problems, as the dissipation becomes negative in certain parts of the domain. Hence, \textit{despite global energy decay, local energy creation cannot be precluded}.
  
\begin{figure}[h!]
    \centering
    \begin{subfigure}[b]{0.45\textwidth}
        \includegraphics[scale=0.2]{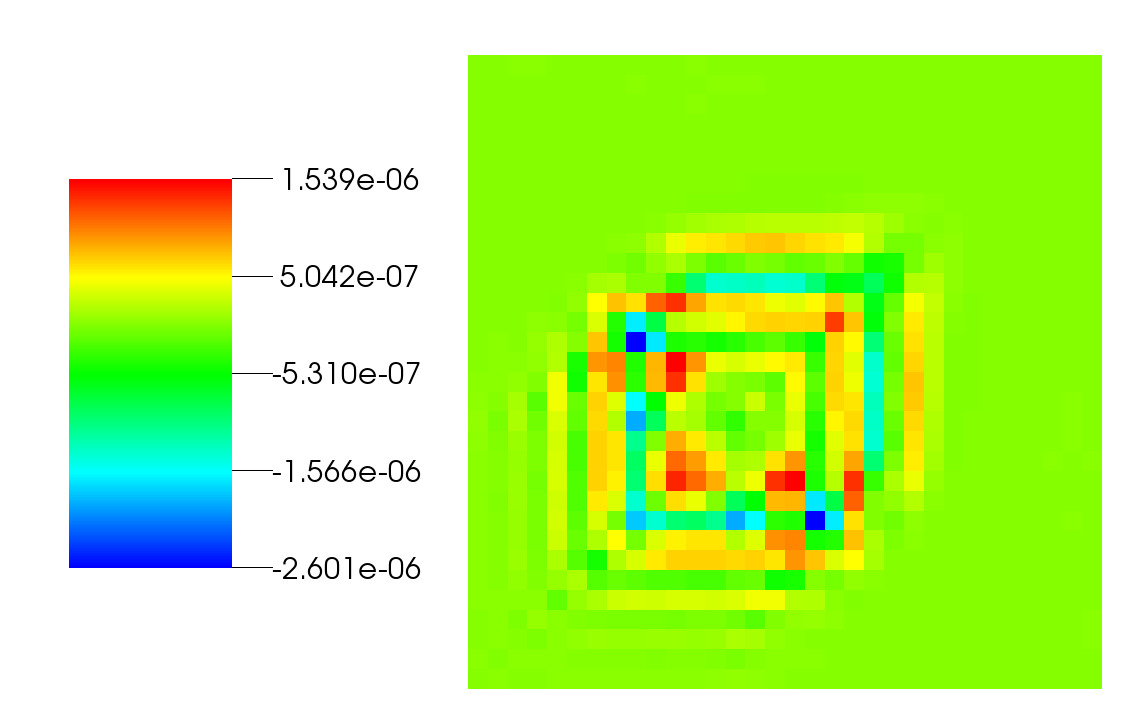}
        \caption{SUPGS, large-scale energy}
    \end{subfigure}
        \begin{subfigure}[b]{0.45\textwidth}
        \includegraphics[scale=0.2]{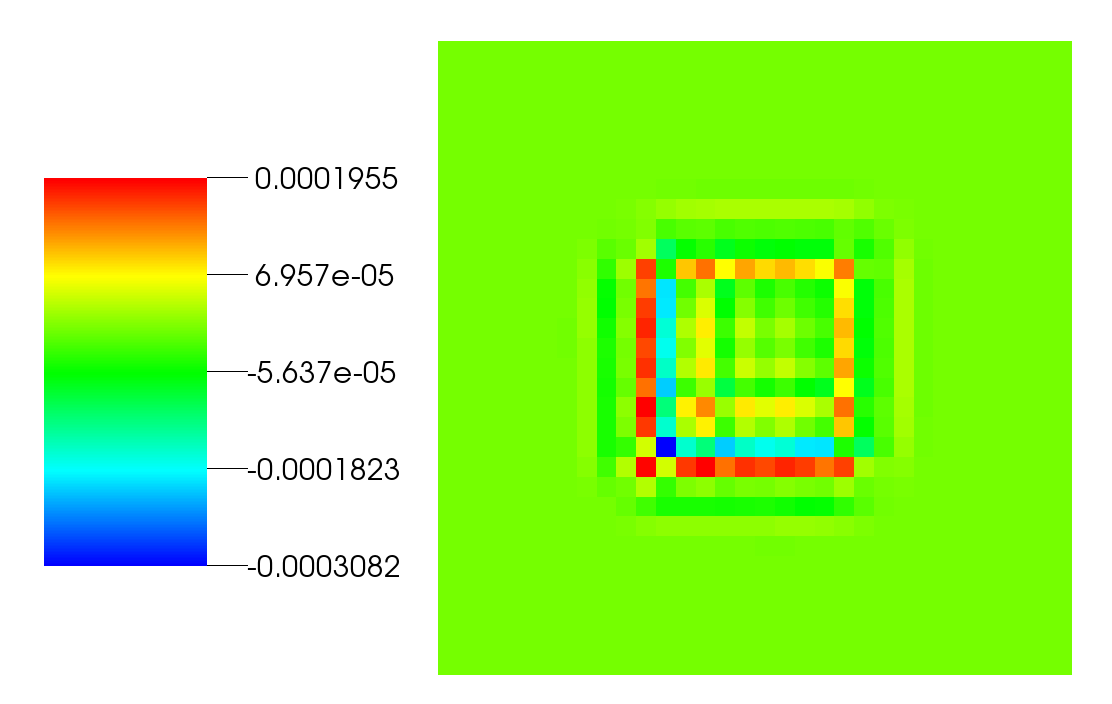}
        \caption{SUPGS, total energy}
    \end{subfigure}
        \begin{subfigure}[b]{0.45\textwidth}
        \includegraphics[scale=0.2]{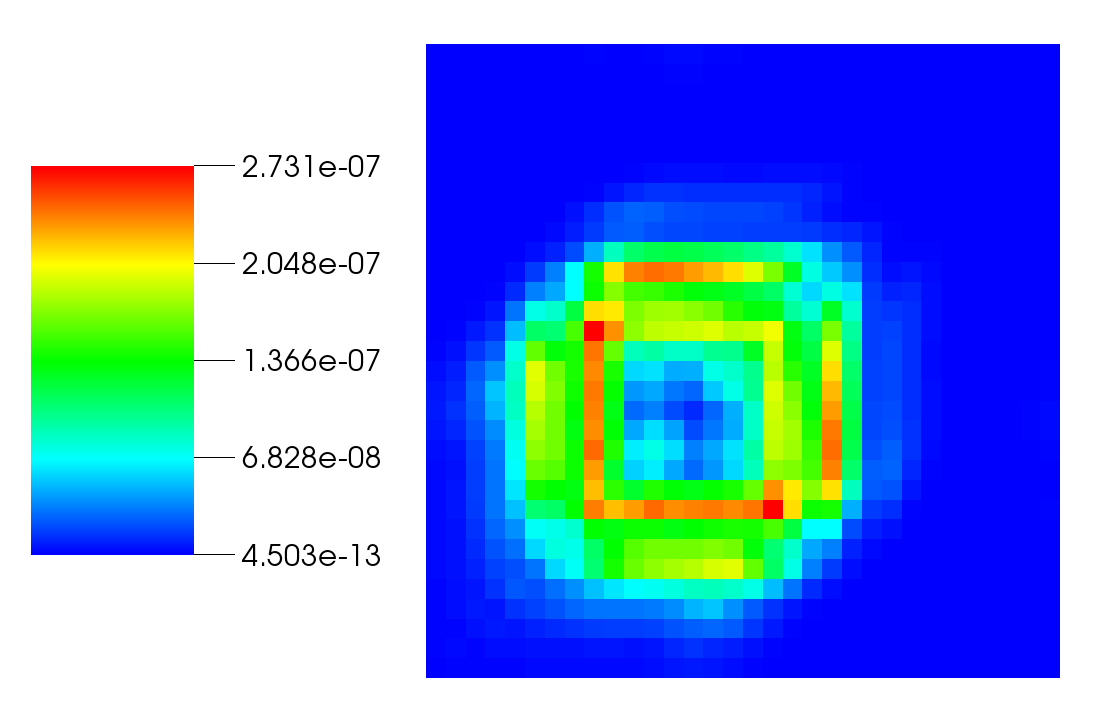}
        \caption{GLSD, total energy}
    \end{subfigure}
    \begin{subfigure}[b]{0.45\textwidth}
        \includegraphics[scale=0.2]{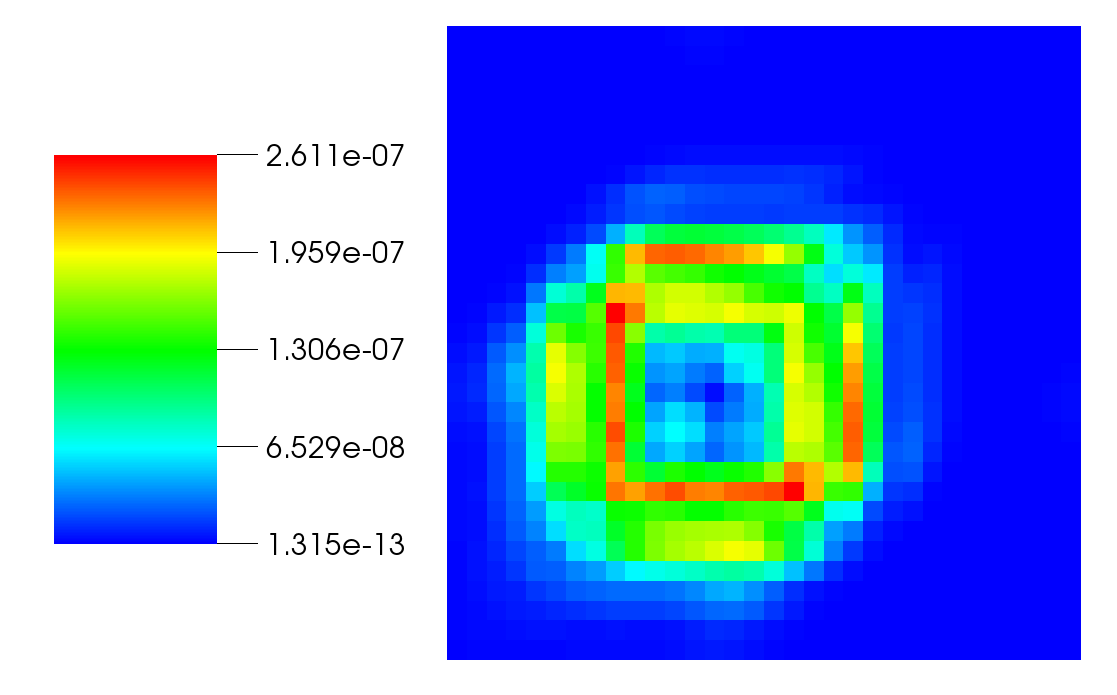}
        \caption{DO, total energy}
    \end{subfigure}
    \caption{Small-scale contribution to energy dissipation on a local scale (at $t=1.0$): (a) the large-scale energy for SUPGS and the total energy for (b) SUPGS, (c) GLSD and (d) DO.}\label{fig:Small-scale dissipation local}
\end{figure}
In the following we further analyze the energy dissipation by considering the contribution of (i) the temporal terms (the last terms on the right-hand side of (\ref{eq:SUPG evolution})) and (ii) the orthogonality term $(\kappa \Delta \phi^h,\phi')$.

\clearpage
\subsection{Temporal-term}
Figures \ref{fig:Temporal contribution to energy dissipation global} and \ref{fig:Temporal contribution to energy dissipation local} show the magnitude of the temporal terms for SUPGS on both global and local level, respectively.  The temporal term of the total energy has larger values than that of the global one. Both energy definitions show negative dissipation, globally as well as locally. Hence contributions of these terms are undesirable.
Comparing with Figure \ref{fig:Small-scale dissipation local} we observe that the temporal has a major contribution to the small-scale dissipation in this case.

\begin{figure}[h!]
    \begin{center}
        \includegraphics[scale=0.65]{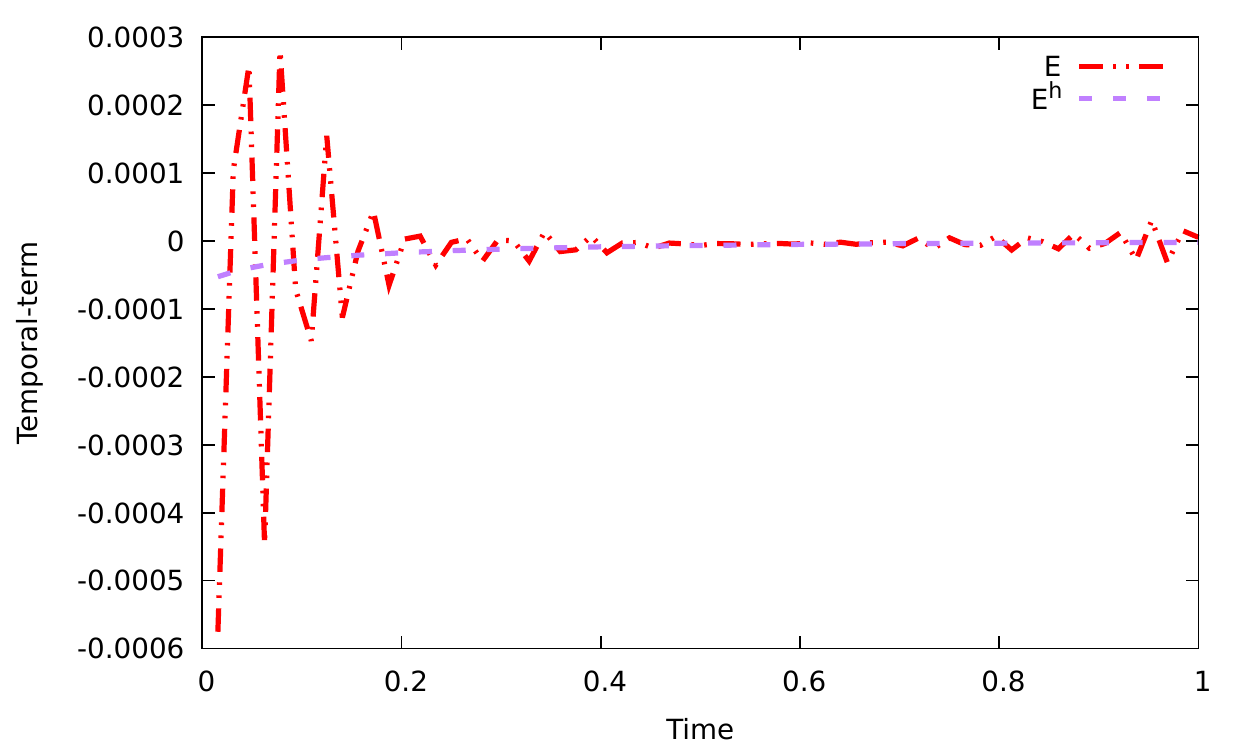}
    \caption{Temporal contribution to small-scale energy dissipation on a global scale for the SUPGS method. The contributions to both the large-scale ($E^h$) and total energy ($E$) are displayed.}\label{fig:Temporal contribution to energy dissipation global}
    \end{center}
\end{figure}

\begin{figure}[h!]
    \begin{center}
    \begin{subfigure}[b]{0.45\textwidth}
        \includegraphics[scale=0.2]{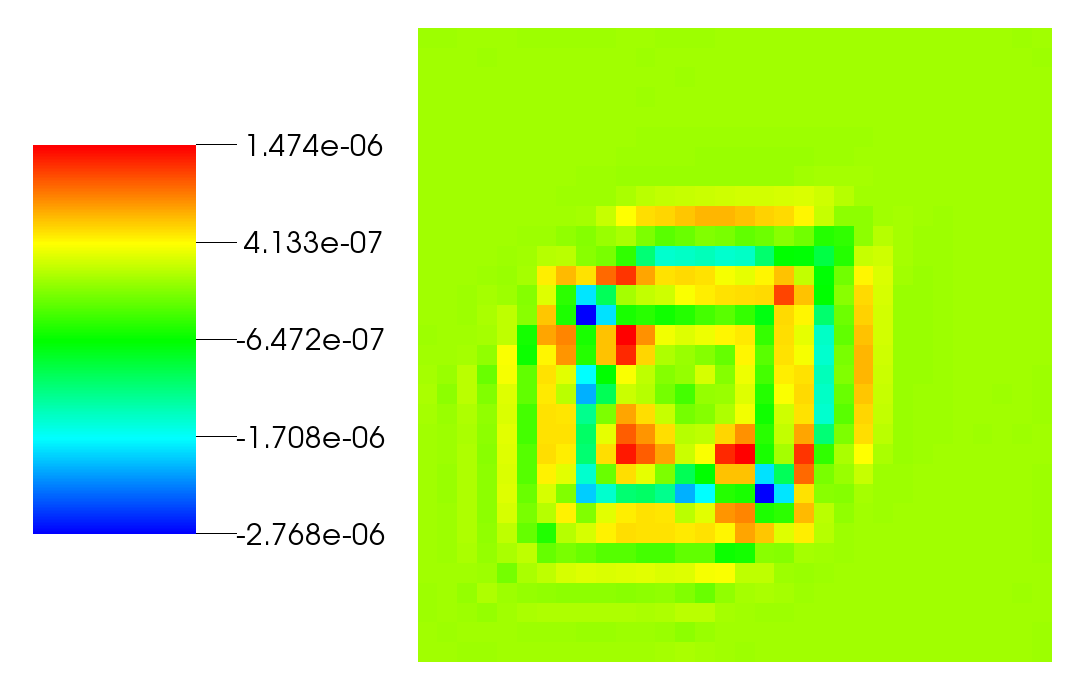}
        \caption{Large-scale energy}
    \end{subfigure}
    \begin{subfigure}[b]{0.45\textwidth}
        \includegraphics[scale=0.2]{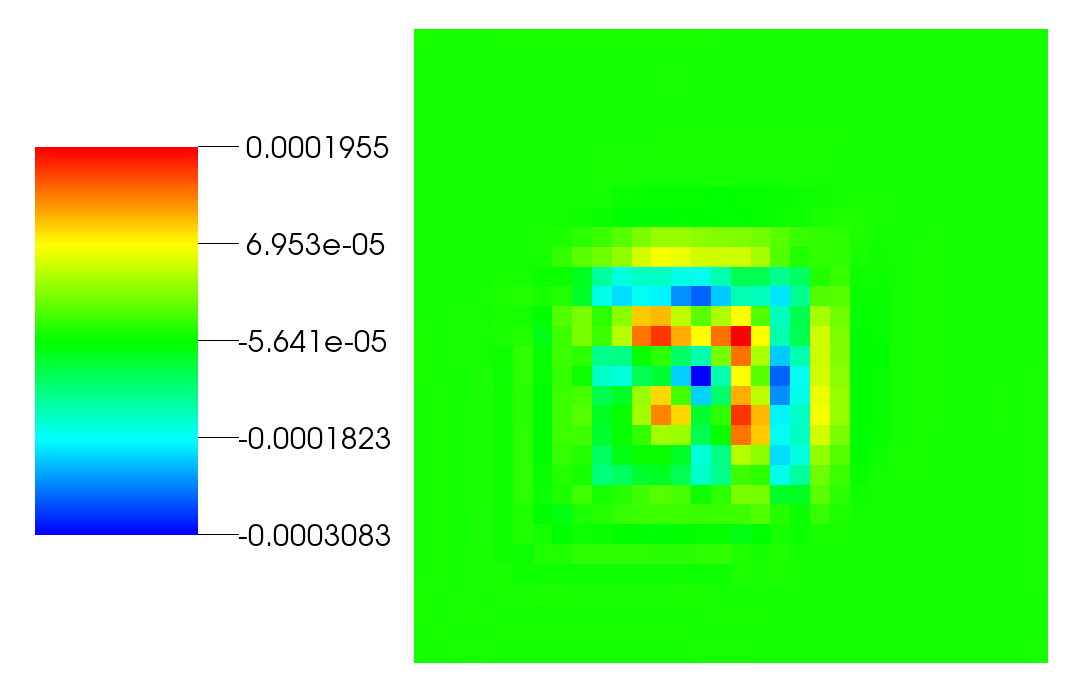}
        \caption{Total energy}
    \end{subfigure}
    \caption{Temporal contribution to small-scale energy dissipation on a local scale for the SUPGS method (at $t=1.0$). The contributions to both the large-scale and total energy are displayed.}\label{fig:Temporal contribution to energy dissipation local}
    \end{center}
\end{figure}

\clearpage

\subsection{Orthogonality-term}
Before continuing we would like to stress that the orthogonality term $\left(\kappa \Delta \phi^h,\phi'\right)$ plays different roles in the formulations.
In case of SUPGS it is  directly an error in the energy behavior, while for GLSD this is an error in the assumed scale separation projector that leads to the correct behavior.
Obviously, for DO the orthogonality term should vanish.

The global and local behavior of the orthogonality term  is displayed in the Figures \ref{fig:ortho term global} and \ref{fig:ortho term local} respectively.
These confirm that the orthogonality term vanishes for the DO formulation.
For the other methods this is not the case. The global orthogonality has an undetermined sign. Moreover, locally the contribution can be negative while the overall contribution is positive.

\begin{figure}[h!]
\begin{center}
        \includegraphics[scale=0.70]{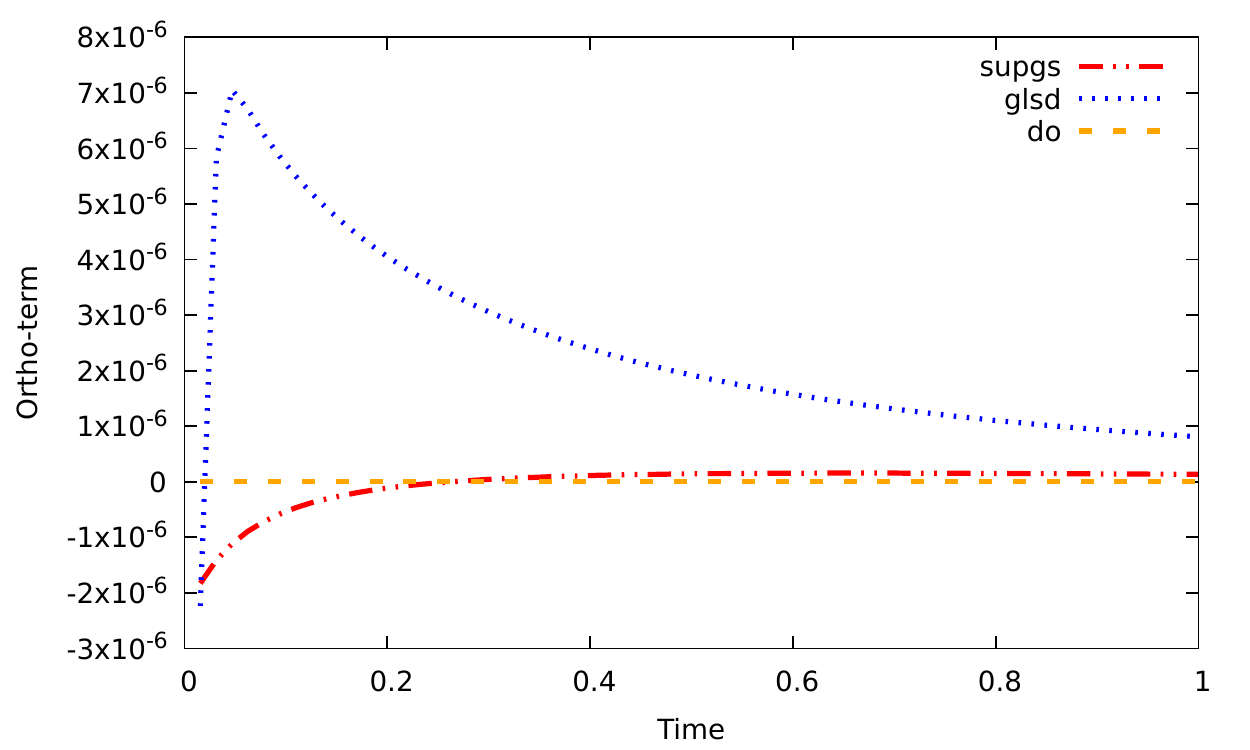}
        \caption{Time evolution of the global orthogonality-term $(\kappa \Delta \phi^h,\phi')_{\tilde{\Omega}}$ for SUPGS, GLSD and DO.}\label{fig:ortho term global}
\end{center}\end{figure}

\begin{figure}[h!]
  \begin{center}
    \begin{subfigure}[b]{0.3\textwidth}
        \includegraphics[scale=0.16]{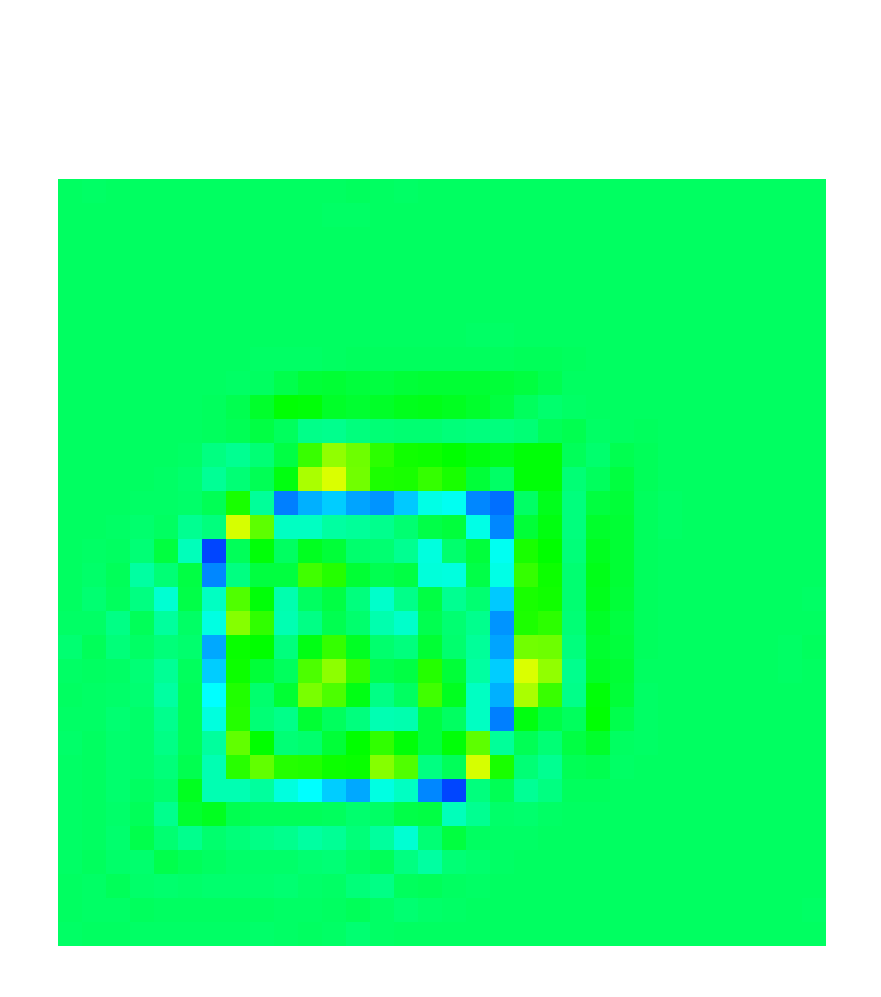}
        \caption{SUPGS}
    \end{subfigure}
    \quad   
    \begin{subfigure}[b]{0.3\textwidth}
        \includegraphics[scale=0.16]{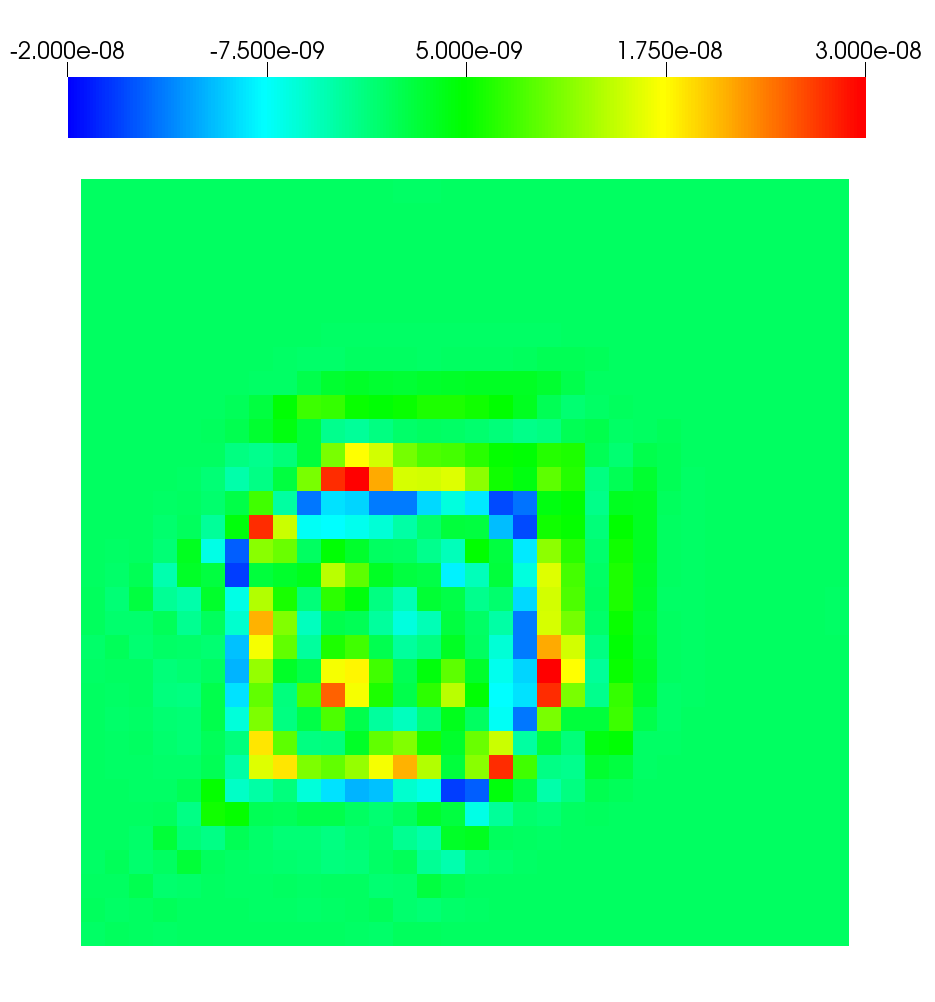}
        \caption{GLSD}
    \end{subfigure}
    \quad  
    \begin{subfigure}[b]{0.3\textwidth}
        \includegraphics[scale=0.16]{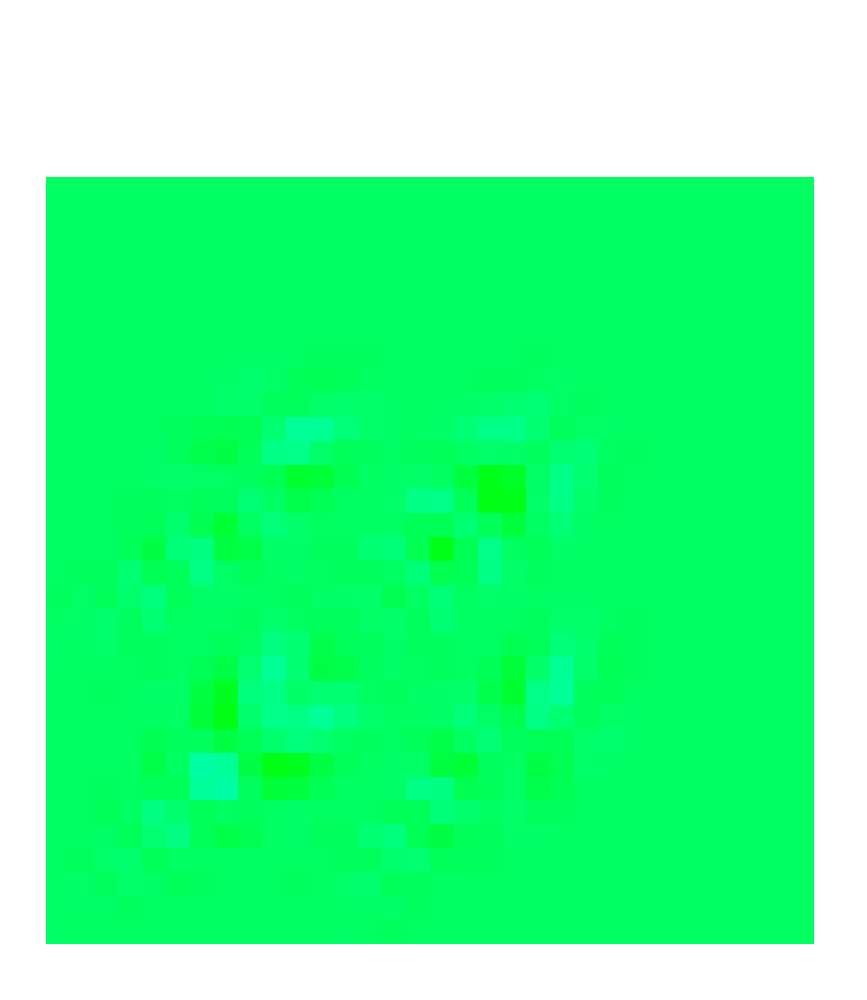}
        \caption{DO}
    \end{subfigure}
    \caption{Local behavior of the orthogonality-term $(\kappa \Delta \phi^h,\phi')_{\tilde{\omega}}$ (at $t=1.0$) for SUPGS, GLSD and DO.}\label{fig:ortho term local}
  \end{center}
\end{figure}
\newpage

\section{Conclusions}\label{sec:conc}
In this work we have proposed an approach to rectify the incorrect-energy behavior of the standard stabilized methods. To this purpose we have employed the concepts of orthogonal small-scales and the dynamic behavior of the small-scales.

This paper takes a road through the various standard weak formulations. The standard Galerkin shows correct-energy evolution but suffers from stability issues. Standard stabilized methods display the opposite. Starting from the variational multiscale approach, we have formulated a design condition to step-by-step remedy the incorrect-energy behavior. The first part towards rectification employs dynamic behavior of the small-scales and henceforth leads to a variational multiscale approach with dynamic small-scales. Next, an orthogonality demand of the large- and small-scales, which can be understood as a $H_0^1$-projection operator, appears. This leads to several options for the variational formulation. It links the form to, both employing dynamic small-scales, the streamline-upwind Petrov-Galerkin method or the Galerkin/least-squares method of which the latter one, in contrast to the former one, displays the energy behavior aimed at. Explicitly enforcing the orthogonality of the large- and small-scales returns us to the variational multiscale framework with the correct-energy behavior. 

The proposed variational formulations which depict correct-energy behavior are:
\begin{itemize}
 \item the Galerkin/least-squares formulation with dynamic small-scales (GLSD)
 \item the approach with dynamic orthogonal small-scales (DO)
\end{itemize}

Numerical results show that the energy convergence of the novel methods displays very similar performance in comparison with the existing stabilized finite element methods. However, the standard methods display both positive and negative small-scale contributions to energy dissipation. The GLSD and the DO method do not suffer from these deficiencies. Furthermore, the numerical results show activity of the unwanted terms in the standard stabilized forms and confirm the enforced orthogonality of the large- and small-scales. The numerical computations have been performed with isogeometric analysis, which is required for a solenoidal velocity field and seems a natural choice when employing orthogonal small-scales. 

This paper serves as an important first step for generalizations in other contexts. Future work will entail a similar methodology for the incompressible Navier-Stokes equations. In particular, we are interested in turbulence computations and applications in science and marine engineering.

\section*{Acknowledgement}\label{sec:ack}

The authors are grateful to Delft University of Technology for its support.

\bibliographystyle{unsrt}
\bibliography{references}

\end{document}